%%%%%%%%%%%%%%
%
%  final version, Oct. 31, 2013
%  sent to the Amer. J. of Math.
%
%
%%%%%%%%%%%%%%%%%%%%%%%%%%%

\magnification=\magstephalf
\input amstex
\documentstyle{amsppt}
\pageheight{9truein}
\pagewidth{6.5truein}
\NoBlackBoxes
\TagsOnRight

\define \fa {\frak A}
\define \fb {\frak B}
\define \fc {\frak C}
\define \fd {\frak D}
\define \fe {\frak E}
\define \fI {\frak I}
\define \fp {\frak P}
\define \fw {\frak W}
\define \fo {\frak O}
\define \fcchi {\fc _{\chi}}
\define \chic {\chi _{\fc}}
\define \chicm {\chi _{\fc ^-}}
\define \Phic {\Phi _{\fc}}

\define \order {\text{\rm ord}}
\define \ord {\text{\rm ord}_v}
\define \ordo {\text{\rm ord}_{v_0}}
\define \ordw {\text{\rm ord}_w}
\define \divi {\text{\rm div}}
\define \Div {\text{\rm Div}}
\define \supp {\text{\rm Supp}}
\define \tr {\text{\rm Tr}}
\define \La {\Lambda}
\define \cS {\Cal S}
\define \difff {\text{\rm Diff}_K}
\define \diff {{\text{\rm Disc}}_K}

\define \ba {\Bbb A}
\define \br {\Bbb R}
\define \bc {\Bbb C}
\define \be {\Bbb E}
\define \bq {\Bbb Q}
\define \bp {\Bbb P}
\define \bfF {\Bbb F}
\define \fq {\bfF _q}
\define \bfp {\bfF _p}
\define \bz {\Bbb Z} 

\define \ka {K_{\ba}}
\define \foba {\fo _{\ba}}
\define \foc {\foba (\fc )}
\define \pia {\pi _{\ba}}
\define \gl {\text{\rm GL}}
\define \gln {\text{\rm GL}_n}
\define \glnp {\text{\rm GL}_{n+1}}
\define \glnm {\text{\rm GL}_{n-1}}
\define \gld {\text{\rm GL}_d}
\define \glnd {\text{\rm GL}_{n-d}}
\define \glnka {\gln (\ka )}
\define \glnpka {\glnp (\ka )}
\define \glnmka {\glnm (\ka )}
\define \gldka {\gld (\ka )}
\define \glnk {\gln (K)}
\define \zetak {\zeta _K}

\define \ux {\bold x}
\define \us {\bold s}
\define \ue {\bold e}
\define \uv {\bold v}
\define \uu {\bold u}
\define \ua {\bold a}
\define \ub {\bold b}
\define \uo {\bold 0}
\define \uy {\bold y}
\define \uY {\bold Y}
\define \uX {\bold X}
\define \uV {\bold V}
\define \ur {\bold r}
\define \uw {\bold w}

\topmatter
\title {On sums involving products and quotients of L-functions over function fields}
\endtitle
\rightheadtext{sums involving L-functions}
\author Jeffrey Lin Thunder
\endauthor
\address Mathematical Sciences Dept., Northern Illinois Univ., DeKalb, IL 
60115
\endaddress
\email jthunder\@ math.niu.edu\endemail

\abstract
We estimate the sum of products or quotients of $L$-functions, where the
sum is taken over all quadratic extensions of given genus over 
a fixed global function field.  Our
estimate for the sum of the quotient of two $L$-functions is analogous to
a result of Schmidt where he estimates the sum of the quotient of
two $L$-series, where the sum is over quadratic extensions  of $\bq$ with
absolute value of the discriminant less than a given bound.
\endabstract

\endtopmatter
\document
\baselineskip=20pt

\head Introduction\endhead

For well over a century zeta functions and $L$-series have played
a prominent role in number theory and are objects of
intense interest in their own right. We mention here just one
example of how they occur. 
A basic notion in Diophantine geometry
is the height, and counting points of given height or height less
than a given bound on various varieties is a subject of much interest.
Early on it was recognized that zeta functions appear
in asymptotic estimates for the number of points of bounded height
in projective spaces. 
A specific instance here is when one counts points of height no
greater than a bound $B$ in $\bp ^2(\overline{\bq})$
that generate a quadratic extension of $\bq$. Schmidt in [9] gave an
asymptotic estimate (in the bound $B$) for the number of such points,
and showed that
the coefficient of the main term in his asymptotic estimate involved
a sum of quotients of $L$-series of the form
$$\sum\Sb [F:\bq ]=2\\ |\text{Disc}(F)|\le m\endSb {L_F(1)\over L_F(3)},$$
where $\text{Disc}(F)$ is the discriminant of $F$ over $\bq$.
In the appendix of [9] he proved asymptotic (in
the parameter $m$) estimates for sums of the form
$$\sum\Sb [F:\bq ]=2\\ |\text{Disc}(F)|\le m\endSb {L_F(s)\over L_F(t)},$$
with certain restrictions on $s,t\in\bc$, of course.
As Schmidt remarks, sums like the one above have received interest before, going
back all the way to Gauss, in fact.

Now for quite some time number theorists have realized that it is
fruitful to also study function fields; says Weil in the foreword
to [13] ``... it goes without saying that the function fields over
finite fields must be granted a fully simultaneous treatment with
number fields ..."  In particular, one can study heights over
any global function field and hopefully gain insight into number fields
in the process. It is in this spirit that the author together with
Widmer in [12] derive asymptotic (in the height)
estimates for the number of points of given height that 
generate an extension
of given degree over a fixed function field, though not in all situations where
we believe such estimates should hold. Together with Kettlestrings, in [7] we
derive such estimates in the remaining previously unproven cases when one
counts points that generate a quadratic extension of a fixed function field.
This is analogous to Schmidt's work cited above, and unsurprisingly requires 
estimates for sums of quotients of $L$-functions over function fields analogous 
to the sums above considered by Schmidt.

Here we will develop machinery (specifically, Propositions 5-7 in \S 4 below)
that allows one to estimate certain sums of
products and quotients of $L$-functions. Loosely speaking, our methods
follow those of Siegel in [10]; they are purposely less ``geometric" and
more ``number-theoretic" in nature. As an example of the use of Propositions 5-7,
we will prove estimates for the sum over the quotient and product of
two $L$-functions. This includes the case of interest above  in [7]
that arises in the counting of points of given height generating a quadratic
extension. Before stating these estimates, we first set some notation.

For a prime $p$ let $\bfp$ denote the finite field with $p$ elements and
let $X$ be transcendental over $\bfp$, so that $\bfp (X)$ is a field
of rational functions. Fix algebraic closures $\overline{\bfp }$ of $\bfp$
and $\overline{\bfp (X)}\supset\overline{\bfp}$ of $\bfp (X)$.
In what follows, by {\it global function field} (or simply {\it function
field}) we mean a finite algebraic extension $K\supseteq\bfp (X)$ contained
in $\overline{\bfp (X)}$. For such a field $K$ we have 
$K\cap
\overline{\bfp}=\bfF _{q_K}$ for some finite field $\bfF _{q_K}$ with
$q_K$ elements; this field is called the field of constants of $K$.
We write $g_K$ for the genus of $K$ and $\zetak$ for the zeta function of
$K$ (defined explicitly below). The $L$-function $L_K$ is given by
$$L_K(q_K^{-s})=(1-q_K^{-s})(1-q_K^{1-s})\zetak (s)={\zetak (s)\over
\zeta _{\bfF _{q_K}(X)}(s)}.$$
It is well-known that $L_K$ is a polynomial of degree $2g_K$ in $q_K^{-s}$ 
and all
its zeros have $\Re (s)=1/2$ (see [11, Chap. V], for example). 
We denote the set of places of $K$ by $M(K)$. 

\proclaim{Theorem 1} Suppose $K$ is a global function field with
field of constants $\fq$,
$m$ is a positive integer and $\epsilon >0$. Suppose $s,t\in\bc$
satisfy i) 
$\Re (s),\ \Re (t)> 3/4+\epsilon$ and $\Re (s)+\Re (t)> 2+2\epsilon$
if $q$ is odd; ii)
$\Re (s),\  \Re (t)> 1/2+\epsilon$ and $\Re (s)+\Re (t)> 3/2+2\epsilon$
if $q$ is even. Then
$$\multline
\sum\Sb [F:K]=2\\ g_F=m,\ q_F=q\endSb 
L_F(q^{-s})L_F(q^{-t})= {2J_Kq^{3-5g_K}\zetak (2s)\zetak (2t)L_K(q^{-s})L_K(q^{-t})
\over q-1}\sigma _1(s,t)
q^{2m}\\
+\cases 
O\big ( 
q^m(1+q^{2m(5/4+\epsilon -\Re (s))})(1+q^{2m(5/4+\epsilon -\Re (t))})\big )
&\text{if $q$ is odd,}\\
O\big ( 
q^m(1+q^{2m(1+\epsilon -\Re (s))})(1+q^{2m(1+\epsilon -\Re (t))})\big )
&\text{if $q$ is even,}\endcases\endmultline$$
where
$$\aligned
\sigma _1(s,t)&=\prod _{v\in M(K)}\Big (1-q^{-2\deg (v)}-(q^{-\deg (v)}-q^{-2\deg (v)})
(q^{-2s\deg (v)}-q^{-2(s+t)\deg (v)}+q^{-2t\deg (v)})\\
&\qquad +(1-q^{-\deg (v)})q^{-(s+t)\deg (v)}\Big )\endaligned$$
and the implicit constants depend only on $K$ and $\epsilon$.
\endproclaim

\proclaim{Theorem 2} Suppose $K$ is a global function field with
field of constants $\fq$,
$m$ is a positive integer and $\epsilon >0$. Suppose $s,t\in\bc$
satisfy i) 
$\Re (s)> 3/4+\epsilon$, $\Re (t)> 1 +\epsilon$ and $\Re (s)+\Re (t)> 2+
2\epsilon$ if $q$ is odd; ii)
$\Re (s)> 1/2+\epsilon$ and $\Re (t)> 1+\epsilon$ if $q$ is even. Then
$$\multline
\sum \Sb [F:K]=2\\ g_F=m,\  q_F=q\endSb {L_F(q^{-s})\over L_F(q^{-t})}=
{2J_Kq^{3-5g_K}\zetak (2s)L_K(q^{-s})\over L_K(q^{-t})(q-1)}\sigma _2(s,t)q^{2m}\\
+\cases O\Big ( 
q^m(1+q^{2m(5/4+\epsilon -\Re (s))}+q^{2m(5/2+2\epsilon -2\Re (t))}+
q^{2m(5/2+2\epsilon -\Re (s)-\Re (t))}\Big )&
\text{if $q$ is odd,}\\
O\big ( q^m(1+q^{2m(1+\epsilon-\Re (s))})\big )
&\text{if $q$ is even,}\endcases\endmultline$$
where
$$\sigma _2(s,t)=
\prod\Sb v\in M(K)\endSb \Big (
1-q^{-2\deg (v)}+q^{-2(s+1)\deg (v)}-q^{-(2s+1)\deg (v)}-q^{-(t+s)\deg (v)}
+q^{-(t+s+1)\deg (v)}\Big )$$
and the implicit constants depend only on $K$ and $\epsilon$.
\endproclaim

\proclaim{Theorem 3} Suppose $K$ is a global function field with
field of constants $\fq$,
$m$ is a positive integer and $\epsilon >0$. Suppose $s,t\in\bc$
satisfy $\Re (s),\  \Re (t)> 1+\epsilon$. Then
$$\multline
\sum \Sb [F:K]=2\\ g_F=m,\  q_F=q\endSb {1\over L_F(q^{-s}) L_F(q^{-t})}=
{2J_Kq^{3-5g_K}\over L_K(q^{-s})L_K(q^{-t})(q-1)}\sigma _3(s,t)q^{2m}\\
+\cases O\big (
q^m(1+q^{2m(5/2+2\epsilon -2\Re (s))}+q^{2m(5/2+2\epsilon -2\Re (t))})\big )
&\text{if $q$ is odd,}\\
O(q^m)&\text{if $q$ is even,}\endcases\endmultline$$
where
$$\sigma _3(s,t)=\prod _{v\in M(K)}\big ( 1-q^{-2\deg (v)}+q^{-(s+t)\deg (v)}-
q^{-(s+t+1)\deg (v)}\big )$$
and the implicit constants depend only on $K$ and $\epsilon$.
\endproclaim

One can view the case where $K$ is a field of rational functions (at least
for odd characteristic) in Theorem 2 as the analog of Schmidt's result over $\bq$
mentioned above. We note that
letting $t=\Re (t)\rightarrow \infty$ in either Theorem 1 or 2 gives the sum of a single
$L$-function. This can be done more directly using our methods, though
the resulting estimate (i.e., error term) is the same. Letting $s=\Re (s)\rightarrow
\infty$ in either Theorem 2 or Theorem 3 gives the sum of the reciprocal of a 
single $L$-function.

\proclaim{Corollary 1} Suppose $K$ is a global function field with
field of constants $\fq$, $m$ is a positive integer and $\epsilon >0$. 
Then for all
$s\in\bc$ with $\Re (s)> 1/2+\epsilon$ if $q$ is even or $\Re(s)>
3/4+\epsilon$ if $q$ is odd, we have
$$\aligned
\sum \Sb [F: K]=2\\ g_F=m,\  q_F=q\endSb L_F(q^{-s})
&=q^{2m}{2J_Kq^{3-5g_K}\zetak (2s)L_K(q^{-s})\over q-1}\\
&\qquad\times \prod \Sb v\in M(K)\endSb
\Big (1-q^{-2\deg (v)}+q^{-2(s+1)\deg (v)}-q^{-(2s+1)\deg (v)}\Big )\\
&\qquad +\cases O\big (q^m(1+q^{2m(5/4+\epsilon-\Re (s))})\big )&\text{if $q$ is odd,}\\
O\big (q^m(1+q^{2m(1+\epsilon-\Re (s))})\big )&\text{if $q$ is even,}
\endcases\endaligned$$
where the implicit constants depend only on $K$ and $\epsilon$.
\endproclaim

\proclaim{Corollary 2} Suppose $K$ is a global function field with
field of constants $\fq$, $m$ is a positive integer and $\epsilon >0$. 
Then for all
$t\in\bc$ with $\Re (t)> 1+\epsilon$  we have
$$\aligned
\sum \Sb [F: K]=2\\ g_F=m,\  q_F=q\endSb {1\over L_F(q^{-t})}
&=q^{2m}{2J_Kq^{3-5g_K}\over L_K(q^{-t})\zetak (2)(q-1)}\\
&\qquad+\cases O\big ( q^mq^{2m(5/2+2\epsilon -2\Re (t))}\big )&\text{if $q$ is odd,}\\
O(q^m)&\text{if $q$ is even,}\endcases\endaligned$$
where the implicit constants depend only on $K$ and $\epsilon$.
\endproclaim

Results of a similar nature to our Corollary 1, though only in odd characteristic, 
have been given by Hoffstein and Rosen [5] (when $K=\fq (X)$ only)
and Fisher and Friedberg [3], with later refinements by Chinta, Friedberg and
Hoffstein in [1]. These works use metaplectic Eisenstein series
and double Dirichlet series. In contrast, our methods here
are rather pedestrian, using only ``classical" machinery such as Dedekind's different
theorem and Hurwitz's genus formula, though to be fair we
also avail ourselves of results such as the Riemann hypothesis for curves over
a finite field which are obviously still unproven in the number field case.
Sums analogous
to those in Corollary 1 but with $\bq$ in place of the function field were
studied by Goldfeld and Hoffstein [4]. 

We remark that one could use these theorems to estimate the number of
quadratic extensions $F$ of $K$ with given genus. Such an
estimate is actually used in our proofs, in fact, and
is another consequence of our analysis (specifically, Proposition 7).
Such estimates are either not extant, or at least difficult
to locate in the present literature, and improve on the quadratic
case of the more general estimates obtained in [2] (over $\fq (X)$ and
only in odd characteristic) and [6].
We also note that another ingredient
of our proofs is a function field analog of the Polya-Vinogradov inequality
(Theorem 4 in \S 3). This is likely of independent interest as well;
we were unable to locate such an analog in the literature.

We end this introduction with some more notation to be used throughout
the remainder and give an elementary result on the prime divisors occurring in a
given effective divisor, a result that is analogous to well-known estimates in
elementary number theory. Section 1 is a short discussion on $L$-functions,
Euler products and how to use M\"obius inversion to express the inverse of
the Artin $L$-function. With this in hand, we can indicate explicitly the type of
sums that must be estimated, which quantities should dominate, and which should
only contribute to ``error terms."
Section 2 is devoted to a 
detailed study of the differents and discriminants that arise from quadratic extensions
and the quadratic extensions that have a particular discriminant. This is rather
routine in odd characteristic, but decidedly less so in characteristic two.
In section 3
we reinterpret our $L$-functions in terms of generators and characters. We then
study multiplicative characters in this context and prove our analog of the
Polya-Vinogradov inequality. Section 4 is the heart of our work. Here we state
and prove our main estimates: Propositions 5, 6 and 7. Using these
we prove Theorems 1, 2 and 3 in the final section.

We will follow the usual conventions that empty sums are interpreted to be zero
and empty products are interpreted to be one.

For a function field $K$ we denote the divisor group, i.e., the free
abelian group on the set of places $M(K)$, by $\Div (K)$. The degree
map on $\Div (K)$, normalized to have image $\bz$, 
will be denoted $\deg _K$ or simply $\deg$
if the field is understood. We will always use capital script German letters
$\fa ,\fb ,\ldots $ to denote divisors, with the sole exception of the
zero divisor $0$. The support of a divisor $\fa$, that is, the (possibly empty)
set of places $v$ for which $\ord (\fa )\neq 0$, will be denoted by
$\supp (\fa )$. We say an effective divisor $\fa$ (that is, a divisor $\fa \ge 0$)
is {\it square-free} if
$\ord (\fa )=1$ for all $v\in\supp (\fa )$.
When we work in the case of characteristic 2 it will prove useful to separate
out the square-free part of an effective divisor. Towards that end,
for any divisor $\fa \ge 0$ we set
$$\fa _1=\sum\Sb v\in\supp (\fa )\endSb v,\qquad \fa _2=\fa -\fa _2.$$
When two effective divisors $\fa$ and $\fb$ have disjoint support, we
write $(\fa ,\fb )=0$.

Each place $v\in M(K)$ has a corresponding order function $\ord :K\rightarrow
\bz \cup \{\infty\}$, whence an (ultrametric) absolute value; these absolute
values lead one to
the adele ring $\ka$ (see [13] for all of the necessary background here). 
For any divisor $\fa\in\Div (K)$ we have the
Riemann-Roch space
$$L(\fa )=\{\alpha\in K\: \ord (\alpha )\ge -\ord (\fa )\ \text{for all
$v\in M(K)$}\},$$
which is a vector space of finite dimension $l(\fa )$ over $\bfF _{q_K}$.
The Riemann-Roch Theorem states that
$$l(\fa )=\deg _K(\fa )+1-g_K+l(\fw _K -\fa ),$$
where $\fw _K\in\Div (K)$ is any divisor in the canonical class (see either
[13, Chapter VI] or [11, Chapter I]).
In particular, if $\deg _K(\fa )\ge 2g_K-1$, then $l(\fa )=\deg _K(\fa )
+1-g_K$, so that the number of effective divisors $\fa\in\Div (K)$ with
fixed degree $\deg (\fa )=m$ is bounded by the product of $q_K^m$ and a function 
depending only on $K$:
$$\sum\Sb \fa\in\Div (K)\\ \fa \ge 0\\ \deg (\fa )=m\endSb 1\ll q_K^m\tag 0$$
for all non-negative integers $m$, where the implicit constant depends only on $K$.
In addition to the Riemann-Roch spaces, we will often consider the
subsets
$$L'(\fa )=\{\alpha\in L(\fa )\: \ord (\alpha )=-\ord (\fa )\ \text{for
all $v\in\supp (\fa )$}\}.$$

The zeta function is given by
$$\zetak (s)=\sum\Sb\fa\in\Div (K)\\ \fa \ge 0\endSb q_K^{-s\deg _K(\fa )}$$
for $\Re (s)>1$, where 
the sum converges absolutely by (0).
The resulting function can be analytically continued to a function
with simple poles at $s=0,1$ and precisely $2g_K$
zeros, all with real part equal to 1/2. See [11, Chap. V], for example.
The case where $K=\fq (X),$ a
field of rational functions, is particularly simple:
$$\zeta _{\fq (X)}(s)={1\over (1-q^{-s})(1-q^{1-s})}.$$

We will often use ``M\"obius inversion" on certain sums
over divisors. The M\"obius function on effective divisors $\fa\in\Div (K)$
is defined exactly as in the classical case: $\mu (\fa )=0$ unless $\fa$ is
square-free, in which case $\mu (\fa )=1$ if $\#\supp (\fa )$ is even and
$\mu (\fa )=-1$ if $\#\supp (\fa )$ is odd.
The inversion arguments all rely on the formula (see [12, Lemma 1])
$$\sum\Sb 0\le\fc\le \fa\endSb\mu (\fc )=\cases 1&\text{if $\fa =0$,}\\
0&\text{otherwise.}\endcases\tag 1$$
In addition to the M\"obius function on effective divisors, we will
also use the phi function:
$$\phi (\fa )=\prod \Sb v\in\supp (\fa )\endSb \big (q_K^{\ord (\fa )}-
q_K^{\ord (\fa )-1}\big ).$$

\proclaim{Lemma 0} Let $K$ be a function field. For
all effective divisors $\fa\in\Div (K)$ and all $\epsilon >0$ we have
$$\#\supp (\fa )\le\epsilon\deg (\fa )+c(\epsilon )$$
for some constant $c(\epsilon )\ge 0$ depending only on $K$ and $\epsilon$, and
$$\sum\Sb 0\le\fb\le\fa\endSb 1\ll q_K^{\epsilon\deg (\fa )},$$
where the implicit constant depends only on $K$ and $\epsilon$.
\endproclaim

\demo{Proof} For the first inequality, let $n$ be the least positive integer
greater than $1/\epsilon$. The number $d(n)$ of places of degree less 
than $n$ certainly depends only on $K$ and $n$ (whence $\epsilon$). Now if
$\#\supp (\fa )\ge d(n)$, then 
$$\deg (\fa )\ge n\#\supp (\fa )-nd(n),$$
so that
$$\epsilon\deg (\fa )+d(n)\ge\#\supp (\fa ).$$
Since the maximum  of $\#\supp (\fa )-\epsilon\deg (\fa )$
over all effective divisors $\fa$ with $\#\supp (\fa )< d(n)$
depends only on $K$ and $n$ (whence $\epsilon$), the first inequality in Lemma 0
follows.

For the second inequality, we first note that the function
$$\theta (\fa )=\sum\Sb 0\le\fb\le\fa\endSb 1$$
satisfies $\theta (\fa +\fb )=\theta (\fa )\theta (\fb )$ whenever $(\fa ,\fb )=0$.
Also, $\theta (nv)=n+1$ for all non-negative integers $n$ and places $v.$ Now
set $c_1=\max _{x\ge 0}\log _{q_K}(x+1)-\epsilon x/2$ and
$c_2=\max \{ c_1,0\}$.
Clearly $c_2$ depends only on $q_K$ and $\epsilon$, so by the first part of
the lemma we see that 
$$\sum\Sb v\in\supp (\fa )\endSb c_2\le (\epsilon /2)\deg (\fa )+c$$
for all effective divisors $\fa$ and
some $c\ge 0$ depending only on $K$ and $\epsilon$. Thus, for
any effective divisor $\fa$
$$\aligned\log _{q_K}\big ( \theta (\fa )\big )
&=\sum _{v\in\supp (\fa )}\log _{q_K}\big (\ord (\fa )+1\big )\\
&\le \sum\Sb v\in\supp (\fa )\endSb c_2+\epsilon \ord (\fa )/2\\
&\le (\epsilon /2)\deg (\fa )+c+(\epsilon /2)\deg (\fa )\\
&=\epsilon\deg (\fa )+c.\endaligned$$
\enddemo

\head 1. $L$-Functions and Euler Products\endhead

Throughout this section we fix a function field $K$ and we set $q=q_K$.
We start with the Euler product representation of the zeta function for
$K$:
$$\zetak (s)=\prod _{v\in M(K)} (1-q^{-s\deg _K(v)})^{-1}.$$
Suppose now that $F$ is a quadratic extension of $K$ with the same
field of constants. We then may
write
$$\aligned \zeta _F (s)&=\prod _{w\in M(F)} (1-q^{-s\deg _F(w)})^{-1}\\
&=\prod _{v\in M(K)}\prod \Sb w\in M(F)\\ w|v\endSb (1-q^{-s\deg _F(w)})^{-1}\\
&=\prod _{v\in M(K)}\prod \Sb w\in M(F)\\ w|v\endSb 
(1-q^{-sf(w)\deg _K(v)})^{-1},\endaligned$$
where $f(w)$ denotes the residue class degree, as usual. We'll
denote the ramification index by $e(w)$.
Since we have a quadratic extension, we have only three possibilities
to consider: there is only one place $w|v$ and it satisfies $e(w)=2,\ f(w)=1$;
there is only one place $w|v$ and it satisfies $e(w)=1,\ f(w)=2$;
there are two places of $F$ lying above $v$, both with ramification
index and residue class degree equal to 1. In the first case here
we set $\chi (F/v)=0$, in the second we set $\chi (F/v)=-1$, and
in the third case $\chi (F/v)=1$. Thus for any place $v\in M(K)$
we have
$$\prod \Sb w\in M(F)\\ w|v\endSb (1-q^{-s\deg _F(w)})^{-1}=
(1-q^{-s\deg _K(v)})^{-1}(1-\chi (F/v)q^{-s\deg _K(v)})^{-1}.$$
Now extend $\chi$ to all effective divisors $\fa\in\Div (K)$ by setting
$$\chi (F/\fa )=\prod\Sb v\in \supp (\fa )\endSb \big (\chi (F/v)\big )^{\ord (\fa )}.$$
Then
$$\aligned
\zeta _F(s)&=\zetak (s)\prod _{v\in M(K)}(1-\chi (F/v)q^{-s\deg _K(v)})^{-1}\\
&=\zetak (s)\sum \Sb \fa\in\Div (K)\\ \fa \ge 0\endSb
\chi (F/\fa )q^{-s\deg _K(\fa )},\endaligned$$
whence
$$L_F(q^{-s})=L_K(q^{-s})\sum \Sb \fa\in\Div (K)\\ \fa \ge 0\endSb
\chi (F/\fa )q^{-s\deg _K(\fa )}.$$

In view of this, we have the following.

\proclaim{Definition} Suppose $K$ is a function field with field of
constants $\fq$ and $F$ is a quadratic extension of $F$ with the
same field of constants. For all $s\in\bc$, set
$$L^*_F(q^{-s})= {L_F(q^{-s})\over L_K(q^{-s})}=
\prod _{v\in M(K)}\big ( 1-\chi (F/v)q^{-s\deg _K(v)}\big )
^{-1}=\sum \Sb \fa\in\Div (K)\\ \fa\ge 0\endSb \chi (F/\fa )
q^{-s\deg _K(\fa )}.$$
\endproclaim

These $L^*_F$ are examples of Artin $L$-series, made particularly simple here
since we are only considering quadratic extensions. 
Though the definition above is a formal one that doesn't address convergence issues,
in fact  $L_F^*$ is a polynomial of degree $2g_F-2g_K$ in $q^{-s}$.
We express the reciprocal of $L_F^*$ with the help of the M\"obius function. 

\proclaim{Lemma 1} Let $K$ be a function field with field of constants $\fq$ and
suppose $F\supset K$ is a quadratic extension with the same field of constants.
For all $s\in\bc$ we have
$$\big ( L^*_F(q^{-s})\big )^{-1}=\sum \Sb \fb\in\Div (K)\\ \fb \ge 0\endSb
\mu (\fb )\chi (F/\fb )q^{-s\deg _K(\fb )}.$$
\endproclaim

\demo{Proof}
Using (1) we have (with $\fc =\fa +\fb$)
$$\aligned L^*_F(q^{-s})
\sum _{\fb \ge 0}\mu (\fb )\chi (F/\fb )q^{-s\deg _K(\fb )}
&=\sum _{\fa \ge 0}\chi (F/\fa )q^{-s\deg _K(\fa )}\sum _{\fb \ge 0}
\mu (\fb )\chi (F/\fb )q^{-s\deg _K(\fb )}\\
&=\sum _{\fc\ge 0}\chi (F/\fc )q^{-s\deg _K(\fc )}\sum _{0\le \fb\le \fc }
\mu (\fb )\\
&=1.\endaligned$$
\enddemo

Given all the above, we see that to estimate a sum of products or quotients of
$L$-functions, we must estimate sums of the form
$$\sum \Sb [F: K]=2\\ g_F=m,\ q_F=q\endSb \chi (F/\fc )\tag 2$$
for effective divisors $\fc\in \Div (K)$.
Note that $\chi (F/\fc )=1$ when $\fc \in 2\Div (K)$, 
except when some place in the support of $\fc$ ramifies in $F$.
Heuristically, one expects sums of the form (2) where $\fc\not\in 2\Div (K)$ 
to be ``small" in some sense, so that the main contribution to our sums of 
$L$-functions should
come from those $\fc$ where $\fc \in 2\Div (K).$
Our goal is therefore two-fold: evaluate these sums when $\fc\in 2\Div (K)$
and show the rest only contribute an error term of smaller magnitude.
This is precisely what we achieve in Propositions 5-7 below. 

\head 2. The Different and Discriminant\endhead

As in the case of number fields, we turn to the discriminant
as a means of conveniently indexing the quadratic extensions
of a given function field $K$. We first dispense with the inseparable case,
however.
Suppose that $[F\: K]=2$ and $q_F=q_K$. If $2|q_K$, then it is possible that
$F$ is not a separable extension. However, in this case we have
$K=\{ \alpha ^2\: \alpha\in F\}$ and $g_F=g_K$ by [11, Proposition III.9.2].
This completely determines the unique such quadratic extension, so for
our purposes here we will only consider separable extensions from now on.

Both [11, Chapter III] and [13, Chapter VIII] contain all the information
we require concerning the different and discriminant.
Here we simply note that every extension $F$
of a given function field $K$ has a canonically given divisor $\difff (F)\in\Div (F)$
called the different. The discriminant $\diff (F)$ is the ``relative norm" of
the different:
$$\diff (F)=\sum\Sb v\in M(K)\endSb\left (\sum \Sb w\in M(F)\\ w|v\endSb \ordw \big (
\difff (F)\big )\deg _F(w)/\deg _K(v)\right )\cdot v\in\Div (K).$$
One readily verifies via the definitions that
$$\deg _K\big (\diff (F)\big )=\deg _F\big (\difff (F)\big ).\tag 3$$
(This is true for the relative norm of any divisor in $\Div (F)$, of course, 
not just the different.)

\proclaim{Proposition 1} Suppose $K$ is a function field with field of
constants $\fq$ and $F$ is a separable quadratic extension
of $K$ with $q_F=q$. 
Then the different $\difff (F)\in \Div (F)$ is an effective divisor with
$$2(g_F-1)=4(g_K-1)+\deg _F\big (\difff (F)\big )
=4(g_K-1)+\deg _K\big (\diff (F)\big ).$$ 
If $2\nmid q$, then both $\difff (F)$ and $\diff (F)$
are square-free effective divisors of even degree.
If $2|q$, then $\difff (F)\in 2\Div (F)$ and $\diff (F)\in 2\Div (K)$.
For any effective divisor $\fa\in\Div (K)$, we have $\chi (F/\fa )\neq 0$
if and only if $\diff (F)$ and $\fa$ have disjoint support.
\endproclaim

\demo{Proof} The equation involving the genera is simply an instance
of the Hurwitz Genus Formula ([11, Theorem III.4.12]) and (3).
Dedekind's Different Theorem [11, Theorem III.5.1] implies that
$\ordw \big (\difff (F)\big )=1$ for all $w\in\supp \big (\difff (F)\big )$
when $q$ is odd. We then have $\ord \big (\diff (F)\big )=1$
for all $v\in\supp\big (\diff (F)\big )$ since in this case $F$ is
a Kummer extension of $K$ (see [11, Proposition III.7.3], for
example).  In the case where $q$ is even $F$ is an Artin-Schreier
extension of $K$; we have $\difff (F)\in 2\Div (F)$ and $\diff (F)\in 2\Div (K)$
by [11, Proposition
III.7.8]. The last part of the Proposition follows from Dedekind's Theorem
and [11, Proposition III.7.8].
\enddemo

Via Proposition 1, we see that our sums (2) may be realized
as follows:
$$\sum\Sb [F: K]=2\\ g_F=m,\ q_F=q\endSb \chi (F/\fc )=
\sum\Sb\fd\in\Div (K)\\ \deg (\fd )=
2m-4g_K+2\endSb \sum\Sb [F:K]=2\\ q_F=q\\ 
\diff (F)=\fd\endSb\chi (F/\fc ).\tag 2'$$
It remains to determine which divisors $\fd\in\Div (K)$ are of the form
$\fd =\diff (F)$ for
some separable quadratic extension $F$ (beyond that which is stated in
Proposition 1), and how many quadratic extensions $F$ have $\diff (F)=\fd$
for a given effective divisor $\fd\in\Div (K)$. 
We introduce a bit more notation.

\proclaim{Notation} Let $\fd\in\Div (K)$ be a square-free effective divisor
of even degree if $q_K$ is odd, and simply a non-zero
effective divisor if $q_K$ is
even. Set
$S(\fd )$ to be the set of separable quadratic extensions
$F$ of $K$ with $q_F=q_K$ and $\diff (F)=\fd$ if $q_K$ is odd, and
$\diff (F)=2\fd$ if $q_K$ is even.
Set $N(\fd )=\# S(\fd )$.
\endproclaim

We will consider the cases of odd and even characteristic separately. 

Suppose first that $K$ is a function field with $q_K=q$ odd.
Since the characteristic of $K$ is odd, every quadratic extension 
$F$ of $K$ is of the
form $F=K(y)$, where $y^2=\omega\in K^{\times}$ and $\omega\neq\alpha ^2$ 
for all
$\alpha\in K$. We say $\omega$ {\it generates} $F$ in this case. One
readily sees that $\omega$ and $\omega '$ generate the same quadratic
extension if and only if $\omega '=\alpha ^2\omega$ for some
$\alpha\in K^{\times}$. 

\proclaim{Proposition 2} Suppose $K$ is a function field with $q_K=q$ odd.
Define $\theta \: \Div (K)\rightarrow\Div (K)\otimes
\bz /2\bz$ by $\theta (\sum n_v\cdot v)=\sum
\pi (n_v)\cdot v,$ where $\pi \: \bz\rightarrow \bz /2\bz$ is the canonical 
map.  Let $\psi$ be
the multiplication by 2 map on the group of divisor classes of degree 0.
The image under $\theta$ of the subset of square-free effective 
divisors of even degree is a subgroup $G_K$ of $\Div (K)\otimes \bz /2\bz$.
The image under $\theta$ of the set of discriminants $\diff (F)$
of quadratic extensions $F\supset K$ with $q_F=q$ together with the zero
divisor is a subgroup of index $\#\ker (\psi )$ in $G_K$.
Further, $N(\fd )=2\#\ker (\psi )$ whenever $S(\fd )\neq\emptyset$.
\endproclaim

\demo{Proof}
In what follows, $\Div _0 (K)$ denotes the group of divisors of degree 0
and ${\Cal P}_K$ denotes the group of principal divisors, so that
$\Div _0(K)/{\Cal P}_K$ is a group of order $J_K$. Set $G_K=\theta \big (
\Div _0(K)\big )$. 

For a $\fd\in\Div _0(K)$,
write $\fd =\fd '+2\fd '',$ where $\fd '$  is a square-free effective
divisor. Then
$\theta (\fd )=\theta (\fd ')$ and 
$0=\deg (\fd ')+2\deg (\fd '')\equiv \deg (\fd ')
\mod 2$. On the other hand, let $\fd =v_1+\cdots +v_n$ be a square-free 
effective divisor of even degree, say $\deg (\fd )=2m$.
Choose a positive integer $k$ and a place $v\not\in
\supp (\fd )$
with $k\deg (v_1)+m=\deg (v).$ Then
$\fe =(2k+1)v_1+v_2+\cdots +v_n-2v\in\Div _0(K)$ with $\theta (\fe )=\fd$.
This proves the first part of the proposition.

Next, for a quadratic extension $F\supset K$ with $q_F=q$ generated
by $\omega$ we have $\theta\big (\diff (K)\big )=\theta \big (\divi (\omega )
\big ),$ where
$\divi (\omega )$ is the principal ideal associated with $\omega$:
$\divi (\omega )=\sum \ord (\omega )\cdot v$. 
(See [11, Proposition III.7.3], for example.)
Let $\psi \: \Div _0(K)/{\Cal P}_K\rightarrow \Div _0(K)/{\Cal P}_K$ 
be the multiplication by 2 map: $\psi (\fa +{\Cal P}_K)=2\fa +{\Cal P}_K$. 
Since the kernel of $\theta$ restricted to $\Div _0(K)$ is clearly
$2\Div _0(K)$, we see that the index
$$[G_K\: \theta ({\Cal P}_K)]=[\theta \big (\Div _0(K)\big )\:\theta (
{\Cal P}_K)]={\#\Div _0(K)/{\Cal P}_K\over\#\psi \big (\Div _0(K)/{\Cal P}_K
\big )}={J_K\over J_K/\#\ker (\psi )}=\#\ker (\psi ).$$

Set $n=\#\ker (\psi )$, 
let $\fa _1=0,\ldots ,\fa _n\in\Div _0(K)$ be representatives of
$\ker (\psi )$ and let $\omega _1=1,\ldots ,\omega _n\in K^{\times}$
with $\divi (\omega _i)=2\fa _i$ for $i=1,\ldots ,n$. Then $\omega ,\omega '
\in K$ generate quadratic extensions $F$ and $F'$ with $\diff (F)=\diff (F')$
if and only if $\omega '=a\omega _i\alpha ^2\omega$ for some
$a\in\fq ^{\times}$, $1\le i\le n$ and
$\alpha\in K^{\times}$. 
Finally, if $\omega /\omega '\in\fq ^{\times}$ and
$\omega /\omega '$ isn't a square (in $K$), then $F\neq F'$, since otherwise
there is an $x\in F=F'$ with $x^2\in\fq ^{\times}$ a non-square (in $\fq$)
so that $q_F=q^2$, contradicting our hypothesis.
Since the squares in $\fq ^{\times}$ form a subgroup of index 2 in
$\fq ^{\times}$, we see that $N(\fd )=2n$ whenever $S(\fd )\neq
\emptyset$. 
\enddemo

We now turn to the case of characteristic two. Again our goal is to determine
$S(\fd )$ and $N(\fd )$; this is achieved in Proposition 3. As opposed to
the case of odd characteristic where $N(\fd )$ is easily determined but
$S(\fd )$ is somewhat mysterious, we find that $S(\fd )$ is almost always
non-empty (see Lemma 5 below) but determining $N(\fd )$ is non-trivial
in characteristic 2. 
Throughout the remainder of this section $K$ is a fixed function field with 
$q_K=q$ even.

Suppose $\omega\in K$ with $\omega\neq \alpha ^2+\alpha$ for all $\alpha\in K$. 
We then have
a quadratic extension $F=K(y)$ where $y^2+y+\omega =0$; we say $\omega$ 
{\it generates} the extension $F$. All (separable)
quadratic extensions $F$ arise in this manner. Similar to the case of odd
characteristic, one readily verifies that $\omega$ and $\omega '$
generate the same quadratic extension if and only if
$\omega '=\omega+\alpha ^2+\alpha$ for some $\alpha\in K$.

For a quadratic extension $F$ generated by $\omega$,
$\diff (F)=\sum (n_v+1)\cdot v$ has the following
properties (see [11, Proposition III.7.8]):
$$\aligned  i&)\  n_v\ge -1\qquad\text{for all $v\in M(K)$;}\\
ii&)\ n_v\equiv 1\mod 2\qquad\text{for all $v\in M(K)$;}\\
iii&)\ -n_v=\max _{\alpha\in K}\{\ord (\omega +\alpha ^2+\alpha )\}\qquad
\text{if $n_v>0$;}\\
iv&)\ n_v=-1\qquad\text{if $\ord (\omega +\alpha ^2+\alpha )\ge 0$ for some 
$\alpha \in K$.}\endaligned
\tag 4$$
Obviously $n_v>0$ only if $\ord (\omega )<0$, and $\diff (F)\in 2\Div (K)$ as
mentioned in Proposition 1. We have the following simple result.

\proclaim{Lemma 3} Suppose $\omega\in K$. If $\ord (\omega )<0$ and odd for 
any place $v\in M(K)$, then $\omega$ generates a quadratic extension $F$
with $\ord (\diff (F))=
1-\ord (\omega )$.\endproclaim

\demo{Proof} If $\alpha\in K$ with $\ord (\alpha )\ge 0$, then 
$\ord (\alpha ^2+\alpha )\ge 0$
and $\ord (\omega )=\ord (\omega +\alpha ^2+\alpha )$. 
If $\ord (\alpha ) <0,$ then 
$\ord (\alpha ^2+\alpha )=
\ord (\alpha ^2)\in 2\bz$. This implies that $\ord (\omega +\alpha ^2+\alpha )=
\min\{\ord (\omega ),
\ord (\alpha ^2)\}$. The lemma follows from this and (4).\enddemo

\proclaim{Lemma 4} Suppose $\fa ,\fb\in\Div (K)$ are effective
divisors with $(\fa ,\fb )=0$. Then
$$\multline
\#\{\alpha\in L(\fa +\fb )\: \ord (\alpha )=-\ord (\fa )\ \text{for all
$v\in\supp (\fa )$}\}\\
=q^{\deg (\fb )+1-g_K}\phi (\fa )\\
+\sum\Sb 0\le\fc\le\fa\\ \deg (\fa +\fb -\fc )<2g_K-1\endSb
\mu (\fc )\big (q^{l(\fa +\fb -\fc )}-q^{\deg (\fa +\fb -\fc )+1-g_K}\big ).
\endmultline$$
In particular, if $\deg (\fb )\ge 2g_K-1$, then
$$\#\{\alpha\in L(\fa +\fb )\: \ord (\alpha )=-\ord (\fa )\ \text{for all
$v\in\supp (\fa )$}\}=q^{\deg (\fb )+1-g_K}\phi (\fa ).$$
\endproclaim

\demo{Proof} We first express $L(\fa +\fb )$ as a disjoint union: 
$$L(\fa +\fb )=\bigcup\Sb 0\le\fc\le\fa\endSb
\{\alpha\in L(\fa -\fc +\fb )\: \ord (\alpha )=-\ord (\fa -\fc )\ \text{for all
$v\in\supp (\fa -\fc )$}\}.$$
For notational convenience, temporarily set
$$L'(\fa -\fc ,\fb )=
\{\alpha\in L(\fa -\fc +\fb )\: \ord (\alpha )=-\ord (\fa -\fc )\ \text{for all
$v\in\supp (\fa -\fc )$}\},$$
so that
$$q^{l(\fa +\fb )}=\sum\Sb 0\le\fc\le\fa\endSb
\# L'(\fa -\fc ,\fb ).$$
We use this together with (1) to get (with $\fe =\fc +\fd$)
$$\aligned
\sum\Sb 0\le\fc\le\fa\endSb \mu (\fc )q^{l(\fa +\fb -\fc )}&=
\sum\Sb 0\le\fc\le\fa\endSb\mu (\fc )\sum\Sb 0\le\fd\le\fa -\fc\endSb
\# L'(\fa -\fc -\fd ,\fb )\\
&=\sum\Sb 0\le\fe\le\fa\endSb \# L'(\fa -\fe ,\fb )
\sum\Sb 0\le \fc\le\fe\endSb \mu (\fc )\\
&=\# L'(\fa ,\fb )\\
&=\#\{\alpha\in L(\fa +\fb )\:\ord (\alpha )=-\ord (\fa )\ \text{for all
$v\in\supp (\fa )$}\}.\endaligned$$

Now by the Riemann-Roch Theorem we have
$$\aligned \sum\Sb 0\le\fc\le\fa\endSb \mu (\fc )q^{l(\fa +\fb -\fc )}&=
\sum\Sb 0\le\fc\le\fa\endSb \mu (\fc )q^{\deg (\fa +\fb -\fc )+1-g_K}\\
&\qquad
+\sum\Sb 0\le\fc\le\fa\\ \deg (\fa +\fb -\fc )<2g-1\endSb\mu (\fc )
\big (q^{l(\fa +\fb -\fc )}-q^{\deg (\fa +\fb -\fc )+1-g_K}\big ).
\endaligned$$
Since
$$\aligned
\sum\Sb 0\le\fc\le\fa\endSb \mu (\fc )q^{\deg (\fa +\fb -\fc )+1-g_K}
&=q^{deg (\fb )+1-g_K}q^{\deg (\fa )}
\sum\Sb 0\le\fc\le\fa\endSb\mu (\fc )q^{-\deg (\fc )}\\
&=q^{deg (\fb )+1-g_K}q^{\deg (\fa )}\prod \Sb v\in\supp (\fa )\endSb
1-q^{-\deg (v)}\\
&=q^{deg (\fb )+1-g_K}\phi (\fa ),\endaligned$$
the lemma follows.
\enddemo

\proclaim{Lemma 5} Suppose $\fd\in\Div (K)$ is a non-zero effective divisor and
$\omega\in L'(\fd _1+2\fd _2)$.
Then $\omega\neq \alpha ^2+\alpha$ for all $\alpha\in K$ and $\omega$ 
generates a quadratic extension $F$ with
$\diff (F)=2\fd$. Further, $\omega ,\omega '\in L'(\fd _1+2\fd _2)$ both 
generate the same quadratic extension $F$ if and only if
$\omega '=\omega +\beta ^2+\beta$ for some $\beta\in L(\fd _2)$. 
There are precisely
$$\multline {2\# L'(\fd _1+2\fd _2)\over q^{l(\fd _2)}}=
{2\phi (\fd )\over q^{l(\fw _K -\fd _2)}}\\
+{2\over q^{l(\fd _2)}}\sum\Sb 0\le\fc\le\fd _1+\fd _2\\
\deg (\fd _1+2\fd _2-\fc )<2g_K-1\endSb\mu (\fc )\big ( q^{l(\fd _1+2\fd _2-
\fc )}-q^{\deg (\fd _1+2\fd _2-\fc )+1-g_K}\big )\endmultline$$
such extensions $F$, where $\fw _K\in\Div (K)$ is any divisor in
the canonical class. In particular, if $\deg (\fd _2)\ge 2g_K-1$, then there
are precisely $2\phi (\fd )$ such extensions.
\endproclaim

\demo{Proof} Let $v$ be any place in the support of $\fd _1+2\fd _2$.  Then
$\ord (\omega )<0$ and odd. As in the proof of Lemma 1 we see that
$\ord (\omega )\neq\ord (\alpha ^2+\alpha )$ for all $\alpha\in K$, so that 
$\omega \neq \alpha ^2+\alpha $ for all
$\alpha\in K$. By construction $\ord (\omega )<0$ only if $\ord (\omega )$ 
is odd. Thus
$\ord (\diff (F))=1-\ord (\omega )=2\ord (\fd )$ for all places $v$ in the
support of $\fd$ by Lemma 3. Moreover, by construction $\ord (\omega )\ge 0$ 
for all places $v$ not in the support of $\fd$. Thus $\diff (F)=2\fd$ by (4).

Now suppose $\omega ,\omega '\in L'(\fd _1+2\fd _2)$ generate the same 
quadratic extension $F$.
Then as noted above, $\omega '=\omega +\beta ^2+\beta$ for some $\beta\in K$. 
If $\ord (\beta )<0$, then as shown in the proof of Lemma 3,
$\ord (\omega ')=\ord (\omega +\beta ^2+\beta )=\min\{\ord (\omega ),
2\ord (\beta )\}$. 
Since $\ord (\omega ')=\ord (\omega )$ for all places $v$ in the support of 
$\fd$ and
$\ord (\omega '),\ord (\omega )\ge 0$ for all other places, we see 
that $v$ here must
be in the support of $\fd$ and $2\ord (\beta )>\ord (\omega )=
-1-2\ord (\fd _2)$. 
Thus, $\ord (\beta )\ge -\ord (\fd _2)$
for all places $v$ in the support of $\fd$, so that $\beta\in L(\fd _2)$.
Conversely, if $\beta\in L(\fd _2 ),$ then $\ord (\omega +\beta ^2+\beta )=
\ord (\omega )$ for all places $v$ in the support of $\fd$ and 
$\ord (\omega +\beta ^2+\beta )\ge 0$ for all other places, so that 
$\omega '=\omega +\beta ^2+\beta\in L'(\fd _1+2\fd _2)$
and both generate the same quadratic extension.  

Finally, we note that the mapping $\beta\mapsto \beta ^2+\beta$ is
two-to-one on $L(\fd _2)$ since $\beta _1^2+\beta _1=\beta _2^2+\beta _2$
if and only if $(\beta _1+\beta _2)^2+\beta _1+\beta _2=0$, i.e.,
$\beta _1+\beta _2=0$ or 1. Thus, there are ${2\# L'(\fd _1+2\fd _2)\over
q^{l(\fd _2)}}$ extensions $F$ with a generator in $L'(\fd _1+2\fd _2)$.
The remainder of the lemma now follows from Lemma 4, upon noting that
$$\aligned
{q^{1-g_K}\phi (\fd _1+2\fd _2)\over q^{l(\fd _2)}}&=
{q^{\deg (\fd _2)+1-g_K}\phi (\fd _1+\fd _2)\over q^{\deg (\fd _2)+1-g_K+
l(\fw _K-\fd _2)}}\\
&={\phi (\fd )\over q^{l(\fw _K-\fd _2)}}.\endaligned$$
\enddemo

\proclaim{Lemma 6} Suppose $\fd$ is an effective divisor and
$F$ is a quadratic extension of $K$ with $\diff (F)=2\fd$. Suppose $S$ is a 
finite non-empty set of places.
Then there is a generator of $F$ in
$$\{\alpha\in L(\fd _1+2\fd _2 +2\fa )\:\ord (\alpha )=1-2\ord (\fd )\ 
\text{for all
$v\in\supp (\fd )\setminus S$}\}$$
for some effective divisor $\fa$ with $\supp (\fa )\subseteq S$.
\endproclaim

\demo{Proof} For notational convenience set $n_v=\ord (\fd _1+2\fd _2)=
2\ord (\fd )-1$ for
$v\in\supp (\fd )$. Let $\omega$ be a generator of $F$ and set
$$M_1(\omega )=\{v\in\supp (\fd )\: \ord (\omega )\neq -n_v\}.$$
If $M_1(\omega )$ is not empty, then we claim that there is a generator $\omega '$ 
of $F$ with $\# M_1(\omega ')<\# M_1(\omega )$. To see why,
let $v'\in\supp (\fd )$ with $\order _{v'} (\omega )\neq -n_{v'}$. 
By (4) there is
an $\alpha\in K$ with $\order _{v'} (\omega +\alpha ^2+\alpha )=-n_{v'}$. 
By the
Weak Approximation Theorem [11, Theorem I.3.1] there is a $\beta\in K$ with
$\order _{v'}(\beta )=0$ and $\ord (\alpha +\beta )=0$ for all other 
$v\in\supp (\fd )$.
Then $\omega '=\omega +(\alpha +\beta )^2+\alpha +\beta$ is a generator of 
$F$ with
$\order _{v'}(\omega ')=-n_v$ and $\ord (\omega ')\ge\min\{\ord (\omega ),0\}$
for all other places $v\in\supp (\fd )$, so that $\# M_1(\omega ')
<\# M_1(\omega )$.
Repeatedly applying this claim yields a generator $\omega $ of $F$ with
$\ord (\omega )=-n_v$ for all places $v\in\supp (\fd )$.

By Lemma 3, if $\ord (\omega )<0$ and $v\not\in\supp (\fd )$, 
then $\ord (\omega )$ must be even. Set
$$M_2(\omega ):=\{v\in M(K)\:\ord (\omega )<0,\ v\not\in\supp (\fd )\cup S\}.$$
If $M_2(\omega )$ is empty, we are done. If not, we claim there is a
generator $\omega '$ of $F$ with $\# M_2(\omega ')<\# M_2(\omega )$ and
$\ord (\omega ')=\ord (\omega )$ for all $v\in\supp (\fd )\setminus S$.
To see why,
let $v'\in M_2(\omega )$. By (4) there is an $\alpha\in K$
with $\order _{v'} (\omega +\alpha ^2+\alpha )\ge 0$. 
By the Strong Approximation Theorem
[11 Theorem I.6.4], there
is a $\beta\in K$ such that $\ord (\alpha +\beta )=0$ for all places 
$v\in M(K)$ with
$\ord (\alpha )<0$ (note that $v'$ is such a place)
and $\ord (\beta )\ge 0$ for any other places $v$ in the
complement of $S$. 
We set $\omega '=\omega +(\alpha +\beta )^2+\alpha +\beta$.
Then $\ord (\omega ')=\ord (\omega )$ for all places $v$ in the intersection of
$\supp (\fd )$ and the complement of $S$,
$\order _{v'}(\omega ')\ge 0$ and $\ord (\omega ')\ge\min\{\ord (\omega ),0\}$
for all other places $v$ in the complement of $S$. Since the
original $M_2(\omega )$ is necessarily a finite set, repeatedly 
applying this claim yields the desired generator.
\enddemo

\proclaim{Lemma 7} Suppose $\fd$ and $\fa$ are effective divisors
with $(\fd ,\fa )=0$. Suppose  $F$ is a quadratic extension of $K$ and
the set
$$\{\alpha\in L(\fd _1+2\fd _2 +2\fa )\: \ord (\alpha )=1-2\ord (\fd )\ 
\text{for
all $v\in\supp (\fd )$}\}$$
contains a generator $\omega$ of $F.$  Then $\diff (F)=2\fd +2\fb$ for
some effective divisor $\fb\le\fa$, and 
$\omega '\in L(\fd _1+2\fd _2 +2\fa )$ generates $F$ if and only if 
$\omega '=\omega +\alpha ^2+\alpha$ for some
$\alpha\in L(\fd _2 +\fa )$. In particular, there are exactly
$q^{l(\fd _2+\fa )}/2$ generators of $F$ in this set.
\endproclaim

\demo{Proof} Suppose $F$ is generated by an $\omega$ in the given set.
By (4), $\supp (\diff (F))\subseteq\supp (\fd )\cup\supp (\fa )$.
By Lemma 3 and the definitions of $\fd _1$ and $\fd _2,$
$\ord (\diff (F))=2\ord (\fd )$ for all $v\in\supp (\fd )$,  
so that $\diff (F)\supseteq 2\fd$.
Also by (4), $2\ord (\diff (F))\le 2\ord (\fa )$ for all $v\in\supp (\fa )$.
This proves the first part of the lemma.

Set $n_v=2\ord (\fd )-1$
for $v\in\supp (\fd )$.
Suppose $\omega '\in L(\fd _1+2\fd _2 +2\fa )$ generates $F$. Then
$\omega '=\omega +\alpha ^2+\alpha$ for some $\alpha\in K$ and
$\ord (\omega ')\ge -n_v$ for all places $v\in\supp (\fd )$. By (4), this implies that
$\ord (\omega ')=-n_v$ for all places $v\in\supp (\fd )$. 
If $\ord (\alpha )<-\ord (\fd _2)$
for some $v\in\supp (\fd )$, then $\ord (\alpha ^2+\alpha )=2\ord (\alpha )\le 
-2\ord (\fd _2) -2=n_v-1$
so that $\ord (\omega ')<-n_v$. Thus $\ord (\alpha )\ge -\ord (\fd _2)$ for all
places $v\in\supp (\fd )$. Similarly, if $v\in\supp (\fa )$ and 
$\ord (\alpha )<-\ord (\fa )$, then
$\ord (\alpha ^2 +\alpha )=2\ord (\alpha )<-2\ord (\fa )$ and
$\ord (\omega +\alpha ^2+\alpha )=2\ord (\alpha )<-2\ord (\fa )$ since $\ord (\omega )
\ge -2\ord (\fa )$.
Thus $\ord (\alpha )\ge -\ord (\fa )$ for all places $v\in\supp (\fa )$ and 
$\alpha\in L(\fd _2+\fa )$.

Conversely, suppose $\alpha\in L(\fd _2 +\fa )$ and let $\omega '=\omega +\alpha ^2
+\alpha$. Then for all places $v\in
\supp (\fd _2)$
we have $\ord (\alpha ^2+\alpha )\ge -2\ord (\fd _2)>-n_v=\ord (\omega )$ so
that $\ord (\omega ')=\ord (\omega )=-n_v$.
For all places $v\in\supp (\fa )$
we have $\ord (\alpha ^2+\alpha )\ge -2\ord (\fa )$ so that
$\ord (\omega ')\ge -2\ord (\fa )$ and $\omega '\in L(\fd _1+2\fd _2 +2\fa )$. 
\enddemo

\proclaim{Proposition 3} For all divisors $\fd >0$ with $\deg (\fd _2)\ge 
2g_K-1$ we have $N(\fd )=2\phi (\fd )$. 
If $\fd >0$ and $v_0\in M(K)\setminus\supp (\fd )$ with
$\deg (v_0)\ge 2g_k-1$, then $N(\fd +v_0)-2\phi (\fd +v_0)=2\phi (\fd )
-N(\fd )$. If $v_0$ is any place with $\deg (v_0)\ge 2g_K-1$, 
then $N(v_0)=2\phi (v_0)$.
\endproclaim

\demo{Proof} Suppose $\fd >0$. Choose
a place $v_0\not\in\supp (\fd )$ with $\deg (v_0)\ge 2g_K-1$.
By Lemma 6, every $F\in S(\fd )\cup S(\fd +v_0)$ has a generator in
$$\{\alpha\in L(\fd _1+2\fd _2+2n_Fv_0)\:\ord (\alpha )=1-2\ord (\fd )\ \text{
for all $v\in\supp (\fd )$}\}$$
for some integer $n_F\ge 0$. Since $S(\fd )\cup S(\fd +v_0)$ is a finite set,
we may choose a positive integer $n$ such that all
$F\in S(\fd )\cup S(\fd +v_0)$ have a generator in the set
$$S=\{\omega\in L(\fd _1+2\fd _2+2nv_0)\:\ord (\omega )=1-2\ord (\fd )\ \text{
for all $v\in\supp (\fd )$}\}.$$

By Lemma 7, if $F$ is any quadratic extension of $K$ with
a generator in $S$, then $\diff (F)=2\fd +2jv_0$ for some $0\le j\le n$.
Conversely, every element $\omega\in S$ is a generator of a quadratic extension
since $\ord (\omega )<0$ and odd for all $v\in\supp (\fd )\neq\emptyset$.
For all $j=0,\ldots ,n$  set
$$S_j=\{ F\in S(\fd+jv_0)\: \text{$F$ has a generator in $S$}\},$$
and for all $j=1,\ldots ,n$, set
$$S_j'=\{ F\in S(\fd +jv_0)\: \ \text{$F$ has a generator in
$L'\big (\fd _1+v_0+2\fd  _2+(2j-2)v_0\big )$}\}.$$
Note that 
$(\fd +jv_0)_1=\fd _1+v_0$ and $(\fd +jv_0)_2=\fd _2+(j-1)v_0$ for $j\ge 1$.
Set $N_j=\# S_j$ for all $j=0,\ldots ,n$ and $N_j'=\# S_j'$ for $j\ge 1$.
Clearly $S_j'\subseteq S_j$ for all $j=1,\ldots ,n$. Also 
$S_0=S(\fd )$ and $S_1=S(\fd +v_0)$ by construction, so that
$N_0=N(\fd ),$ $N_1=N(\fd +v_0)$ and $N_j'\le N_j$ for $j\ge 1$.

Since $\deg (v_0)\ge 2g_K-1$ and $n>0$ we have 
$$\# S=q^{2n\deg (v_0)+1-g_K}\phi (\fd _1+2\fd _2)=q^{2n\deg (v _0)+\deg (\fd 
_2)+1-g_K}\phi (\fd )$$ 
by Lemma 4,
so that by Lemma 7
$$\aligned
\sum _{j=0}^nN_i={2\# S\over q^{l(\fd _2+nv_0)}}&={2q^{2n\deg (v_0)
+\deg (\fd _2)+1-g_K}\phi (\fd )\over q^{n\deg (v_0)+\deg (\fd _2)+1-g_K}}\\
&=
2q^{n\deg (v_0)}\phi (\fd ).\endaligned\tag 5$$
By Lemma 5, for all $j=2,\ldots ,n$ we have
$$\aligned N_j'&=2\phi (\fd +jv_0)\\
&=2\phi (\fd )\phi (jv_0)\\
&=2\phi (\fd )\big ( q^{j\deg (v_0)}-q^{(j-1)\deg (v_0)}\big ).
\endaligned\tag 6$$
Therefore
$$\aligned
\sum _{j=0}^n N_j\ge N_0+N_1'+\sum _{j=2}^n N_j'&=
N_0+N_1'+\sum _{j=2}^n2\phi (\fd )\big ( q^{j\deg (v_0)}-q^{(j-1)\deg (v_0)}
\big )\\
&=N_0+N_1'+2\phi (\fd )\big ( q^{n\deg (v_0)}-q^{\deg (v_0)}\big ).
\endaligned\tag 7$$

Suppose $\deg (\fd _2)\ge 2g_K-1$. Then since $(\fd +v_0)_2=\fd _2$ we have
(6) for $j=1$ as well. Combining this with (5) and (7) yields
$$\gathered 2\phi (\fd )q^{n\deg (v_0)}\ge N_0+2\phi (\fd )\big (
q^{n\deg (v_0)}-1\big ),\\
2\phi (\fd )\ge N_0.\endgathered$$
But $N_0=N(\fd )$ by construction and $N(\fd )\ge 2\phi (\fd )$ by Lemma 5.
Hence $N_0=N(\fd )=2\phi (\fd )$.

Now suppose simply that $\fd >0$. Since $\deg (\fd _2+(j-1)v_0 )\ge 
(j-1)\deg (v_0)\ge 2g_K-1$ for all $j=2,\ldots ,n$, by what we have 
already proven we have
$$ N_j'=N_j=N (\fd +jv_0)=2\phi (\fd +jv_0)
=2\phi (\fd )\big ( q^{j\deg (v_0)}-q^{(j-1)\deg (v_0)}\big )$$
for all $j=2,\ldots ,n.$
Thus by (5), 
$$\aligned
2q^{n\deg (v_0)}\phi (\fd )&= N_0+N_1+\sum _{j=2}^n N_j\\
&= N(\fd ) +N(\fd +v_0)+2\phi (\fd )\big ( q^{n\deg (v_0)}-q^{\deg (v_0)}
\big ),\endaligned$$
so that
$$\aligned 2\phi (\fd ) -N(\fd )&= N(\fd +v_0)-2\phi (\fd )q^{\deg (v_0)}
+2\phi (\fd )\\
&=N(\fd +v_0)-2\phi (\fd +v_0).\endaligned$$

Finally, if $\fd =0$ we have the same argument as above except that now
not all elements of $S=L(2nv_0)$ need generate a quadratic extension and
we have one quadratic extension with a generator in $\fq=L(0)$ (the 
field obtained from $K$ be extending the field of constants to $q^2$).
This extension has $q^{l(nv_0)}/2$ generators in $S$ by Lemma 7.
Any $\beta\in L(nv_0)$ gives a $\beta ^2+\beta\in L(2nv_0)$
and vice-versa. Once again, $\beta\mapsto \beta ^2+\beta$ is a
two-to-one map on $L(nv_0)$, so there are $q^{l(nv_0)}/2$ elements of
$S$ that do not generate a quadratic extension of $K$. Therefore (5) is 
replaced by
$$2+\sum _{j=1}^nN_j={2\# S\over q^{l(nv_0)}}=2q^{l(2nv_0)-l(nv_0)}
=2q^{n\deg (v_0)}.$$
Using $N_j=2\phi (jv_0)$ for $j\ge 2$ as above yields
$$2+N(v_0)+2(q^{n\deg (v_0)}-q^{\deg (v_0)})=2q^{n\deg (v_0)},$$
whence $N(v_0)=2q^{\deg (v_0)}-2=2\phi (v_0).$
\enddemo

\proclaim{Corollary} If $g_K=0$, then $N(\fd )=2\phi (\fd )$ for 
all $\fd >0$. If $g_K>0$, then $N(\fd )=2\phi (\fd )$ if no place
$v\in\supp (\fd )$ has $\deg (v)<2g_K-1$ or if $\deg (\fd _2)\ge 2g_K-1.$
For all other divisors $\fd >0$ we have
$$|N(\fd )-2\phi (\fd )|\ll 1,$$
where the implicit constant depends only
on $K$.
\endproclaim

\demo{Proof} The case where $g_K=0$ follows immediately from Propositions 1
and 3. Suppose $g_K>0$ and write 
$$\fd '=\sum\Sb v\in\supp (\fd )\\ \ord (\fd )>1\ \text{or}\ \deg (v)<2g_K-1
\endSb \ord (\fd )\cdot v.$$
Note that $\fd _2 =\fd '_2$ by construction. According to Proposition 3, 
$$|N(\fd )-2\phi (\fd )|=|N(\fd ')-2\phi (\fd ')|.$$
Now if $\deg (\fd _2)=\deg (\fd '_2)\ge 2g_K-1$, then $N(\fd ')=2\phi (\fd ')$
by Proposition 3. If $\deg (\fd '_2)<2g_K-1,$ then $\fd '$ is an effective 
divisor with $\deg (\fd ')\ll 1$. The Corollary follows in this case
from (0). 
\enddemo

\head 3. Characters and Generators\endhead

In this section we reinterpret $\chi (F/\fc )$
in terms of generators and certain characters. We then look at certain
character sums and prove a result akin to the Polya-Vinogradov inequality
for incomplete multiplicative character sums.

As before we fix a function field $K$ and write $q_K=q$.
When $F$ is a separable quadratic extension of $K$ and $\fc\in\Div (K)$ is
an effective divisor, we desire a description of $\chi (F/\fc )$ in
terms of a generator of $F$.
Recall that we say a (separable) quadratic extension $F$ of $K$ is 
generated by an $\omega\in K$ if $F$ is the splitting field of the polynomial
$P_{\omega}(Y)$ given by
$$P_{\omega}(Y)=\cases Y^2-\omega&\text{if $q$ is odd,}\\
Y^2+Y+\omega&\text{if $q$ is even.}\endcases$$

For any place $v\in M(K)$, set  $K_v$ to be the topological completion of $K$
at the place $v$ and set $\fo _v=\{\alpha _v\in K_v:\ord (\alpha _v )\ge 0\}.$
Let $\pi _v\in\fo _v$ with $\ord (\pi _v)=1$. Then $\pi _v\fo _v$ is
the unique maximal ideal in the discrete valuation ring $\fo _v$.
The residue class field $\fo _v/\pi _v\fo _v$ is isomorphic to
$\bfF _{q^{\deg (v)}}$ by definition/construction. Now for any $\alpha _v\in
\fo _v$ define
$$ P_{\alpha _v}(Y)=\cases Y^2-\alpha _v&\text{if $q$ is odd,}\\
Y^2+Y+\alpha _v&\text{if $q$ is even.}\endcases$$
Define $\chi _v(\alpha _v)$ to be 0,  1 or -1, depending on
whether $P_{\alpha _v}$ is a square, a product of two distinct 
linear factors, or is irreducible over $\fo _v/\pi _v\fo _v$, respectively.

Recall that we denote the adele ring by $\ka$. 
For a non-zero effective divisor $\fc\in\Div (K)$ set 
$$\foc =\{ (\alpha _v)\in\ka \: \alpha _v\in\fo _v\ 
\text{for all places $v\in\supp (\fc )$}\}$$
and 
$$\pia (\fc )=\{ (\alpha _v)\in\foc \: \alpha _v\in\pi _v^{\ord (\fc )}\fo _v\ 
\text{for all $v\in\supp (\fc )$}\}.$$
 Then $\chi _v(\alpha _v)$ is
defined for all $v\in\supp (\fc )$ whenever the adele $(\alpha _v)\in\foc$, and
we set
$$\chi _{\fc}\big ((\alpha _v)\big )=\prod _{v\in\supp (\fc )}
\chi _v(\alpha _v)^{\ord (\fc )}.$$

\proclaim{Proposition 4} Suppose $K$ is a function field, $\fc\in\Div (K)$ 
with $\fc >0$ and $F$ 
is a quadratic extension
of $K$ satisfying $q_F=q_K$ and $\big (\fc ,\diff (F)\big )=0.$
Suppose $F$ is generated by $\omega\in\foc$. If $q_K$ is even, then
$\chi (F/\fc )=\chi _{\fc}(\omega )$.
If $q_K$ is odd, then $\chi (F/\fc )=\chi _{\fc}(\omega )$ if
$\ord (\omega )=0$ for all $v\in\supp (\fc )$, i.e., if $\omega +\pia (\fc _1)$
is a unit in the quotient ring $\foc /\pia (\fc _1)$.
\endproclaim

\demo{Proof}
Suppose $F$ is a quadratic extension generated by $\omega\in\foc$ and let
$y\in F$ satisfy $y^2=\omega$ if $q$ is odd, or
$y^2+y=\omega$ if $q$ is even. 
The formal derivative of $P_{\omega}$ satisfies
$$P'_{\omega}(y)=\cases 2y&\text{if $q_K$ is odd,}\\
1&\text{if $q_K$ is even.}\endcases$$
If $v\in\supp (\fc )$, then 
no place $w\in M(F)$ lying above $v$ ramifies by Proposition 1, so that
$0=\ordw \big (P'_{\omega}(y)\big )=\ordw \big (\difff (F)\big )$ and
$\{1,y\}$ is an integral basis for $F$ over $K$ at the place $w$ by
[11, Theorem III.5.10]. The proposition now follows from 
[11, Theorem III.3.7].
\enddemo

With this in hand, we now make some pertinent observations about
$\chic$ and $\foc$.

\proclaim{Lemma 8} Suppose $K$ is a function field and $\fc\in\Div (K)$
with $\fc >0$. Then
$$\foc /\pia (\fc )\cong\prod _{v\in\supp (\fc )}\fo _v /\pi _v^{\ord (\fc )}
\fo _v$$
and $\dim _{\bfF _{q_K}}\big (\foc /\pia (\fc)\big )=\deg (\fc )$.
If $\fa$ is any divisor with $\supp (\fa )\cap\supp (\fc )=\emptyset$, then
$L(\fa )\subset\foc$, and if $\deg (\fa )\ge\deg (\fc )+2g_K-1,$ then
the image of $L(\fa )$ under the canonical map $\foc\mapsto\foc /\pia (\fc )$
is $\foc /\pia (\fc )$.
\endproclaim

\demo{Proof} The first statement is obvious and the second follows immediately
since $\dim _{\bfF _{q_K}}(\fo _v/\pi _v\fo _v)=\deg (v)$ by definition.    

Let $\theta _{\fc}\:\foc \rightarrow \foc /\pia (\fc )$ denote
the canonical map. Suppose $\fa$ is a divisor with $\supp (\fa )\cap\supp (\fc )
=\emptyset .$ Then $L(\fa )\subset\foc$ directly from the definitions and
$$\aligned L(\fa )\cap\ker (\theta _{\fc})&=L(\fa )\cap \pia (\fc )\\
&= \{\alpha\in L(\fa ):\ord (\alpha )\ge \ord (\fc )\ 
\text{all $v\in\supp (\fc )$}\}\\
&=L(\fa -\fc ).\endaligned$$
Suppose  further that
$\deg (\fa )\ge\deg (\fc )+2g_K-1$. Then the dimension of
the image $\theta _{\fc}\big ( L(\fa )\big )$ is 
$l(\fa )-l(\fa -\fc )=\deg (\fc )$
by the Riemann-Roch Theorem. This completes the proof.
\enddemo

\proclaim{Lemma 9} Suppose $K$ is a function field and
$a \in\bfF _{q_K}$ such that $a\neq b^2$ if $q_K$ is odd or
$a\neq b^2+b$ if $q_K$ is even, for all 
$b\in\bfF _{q_K}$. Then for all places $v\in M(K)$ we have
$$\chi _v(a)=\cases 1&\text{if $\deg (v)$ is even,}\\
-1&\text{if $\deg (v)$ is odd.}\endcases$$
In particular, for any divisor $\fc >0$ we have
$\chi _{\fc}(a)=(-1)^{\deg (\fc )}$ and if $\deg (\fc )$ is odd
$$\sum\Sb [F:K]=2\\ g_F=m,\ q_F=q_K\endSb \chi (F/\fc )=0$$
for all $m\ge 0$.
\endproclaim

\demo{Proof}  Set $q=q_K$.
Let $v\in M(K)$ and choose a $\pi _v\in\fo _v$ as above.
We have $\bfF _{q^{\deg (v)}}\cong \fo _v/\pi _v\fo _v$ and
$\fq\cap\pi _v\fo _v=\{ 0\}$ since $\ord (b)=0\neq 1$ for all
non-zero $b\in\fq$. Now
$\{b +\pi _v\fo _v\: b\in\fq\}\cong \fq$ is a subfield of the
residue class field $\fo _v/\pi _v\fo _v$. As such, $Y^2-b$ (if
$q$ is odd) or
$Y^2+Y+b$ (if $q$ is even)
splits over the residue class field $\fo _v/\pi _v\fo _v$
for $b\in\fq$ exactly when the residue class field contains the splitting
field; this splitting field is (isomorphic to) $\bfF _{q^2}$ when $b=a$.  
On the other hand, this is possible
exactly when $\bfF _{q^2}$ is (isomorphically)
a subfield of $\bfF _{q^{\deg (v)}}$. 
Thus $\chi _v(a)=-1$ for all places $v\in M(K)$
of odd degree and
$\chi _v(a)=1$ for all places $v\in M(K)$ of even degree.

Now suppose $F$ is a (separable) quadratic extension of $K$ with $q_F=q$
and say $F$ is generated by $\omega\in K$. Then $a\omega$ 
generates a quadratic extension $F'$ with $\diff (F)=\diff (F')$ in
odd characteristic and $a+\omega$ generates an $F'$ with $\diff (F)=
\diff (F')$ in characteristic 2. We can't have $F'=F$ in either case, since
then $q_F=q^2$. By what we have shown, if $\deg (\fc )$ is odd, then
$\chi (F'/\fc )=-\chi (F/\fc )$. Finally, $b^2\omega$ in odd characteristic
and $\omega +b^2+b$ in characteristic 2 both generate the original
quadratic extension $F$ for all $b\in\fq$. This shows the final equation in Lemma 9.
\enddemo

We now turn to a more general look at characters. For our purposes here,
a character on an abelian group $G$ will mean a group homomorphism from $G$ into
$\bc ^{\times}$. If $H\subseteq G$ is any finite subgroup and
$\chi$ is a character on $G$, then we have
$$\sum \Sb h\in H\endSb\chi (h)=\cases |H|&\text{if $H\subseteq\ker (\chi )$,}\\
0&\text{otherwise.}\endcases\tag 8$$
The proof is trivial. Let $h_0\in H;$ then
$$\sum\Sb h\in H\endSb \chi (h)=\sum\Sb h\in H\endSb\chi (h+h_0)=\chi (h_0)
\sum\Sb h\in H\endSb\chi (h).$$
Now suppose $R$ is a finite commutative ring with identity.
A {\it multiplicative character} on $R$ is a function $\Phi\:
R\rightarrow \bc$ such that the restriction to the multiplicative
group of units $R^{\times}$ is a group character and
$\Phi =0$ on the complement.
When the restriction
to $R^{\times}$ is the trivial character, we call $\Phi$ 
the principal character on $R$. We say the multiplicative character
$\Phi$ is {\it primitive} if for all non-zero ideals $I\subseteq R$
there are units $u,v\in R^{\times}$ with $u\equiv v\mod I$ and
$\Phi (u)\neq\Phi (v)$.

We readily see that $\chic$ is a character on $\foc$ when $q_K$ is even, with
$\ker (\chic )\supseteq\pia (\fc _1)$.
When $q_K$ is odd, we readily see that
$\Phic \big ((\alpha _v) +\pia (\fc _1)\big )=\chic \big ((\alpha _v )\big )$ 
defines a multiplicative character on the quotient ring $\foc /\pia (\fc _1)$.

In what follows we will take the convolution of a multiplicative 
character $\Phi$ on a finite commutative ring with identity $R$ with
a fixed character $\chi$ on the additive structure of $R$:
$$\tau (\Phi ):=\sum\Sb r\in R\endSb \Phi (r)\chi (r).$$
We require the 
character $\chi$ to be non-trivial on all non-zero principal ideals; we assume
this is the case in what follows. Thus, for all non-zero $h\in R$ we have
$$\sum\Sb g\in R\endSb \chi (gh)=0$$
by (8).

\proclaim{Lemma 10} Suppose $R$ is a finite commutative ring with identity
and $\Phi$ is a multiplicative character on $R$.  Then for all 
$u\in R^{\times}$ we have
$$\sum _{r\in R}\overline{\Phi (r)}\chi (ur)=
\Phi (u)\tau (\overline{\Phi }),$$
where the overline denotes complex conjugation.\endproclaim

\demo{Proof} Setting $r'=ur$ below, we have
$$\aligned \Phi (u)\sum _{r\in R}\Phi (r)\chi (ur)&=
\sum _{r\in R}\Phi (ur)\chi (ur)\\
&=\sum _{r'\in R}\Phi (r')\chi (r')\\
&=\tau (\Phi ).\endaligned$$
The lemma follows upon replacing $\Phi$ above by $\overline{\Phi }$ and
noting that $\overline{\Phi (u)}=\big (\Phi (u)\big )^{-1}$.
\enddemo

\proclaim{Lemma 11} Suppose $R$ is a finite commutative ring with identity and
$\Phi$ is a primitive multiplicative character on $R$. Suppose $I$ is a 
non-zero ideal of $R$ and $r_0\in R$. Then
$$\sum \Sb r\in R\\ r\equiv r_0\mod I\endSb \Phi (r)=0.$$
\endproclaim

\demo{Proof} Let $u,v\in R^\times$ with $u\equiv v\mod I$ 
and $\Phi (u)\neq\Phi (v)$. Then $w=u^{-1}v$ satisfies $\Phi (w)\neq 1$ 
and $w\equiv 1\mod I$. We thus have (with $r'=rw^{-1}$)
$$\aligned
\sum \Sb r\in R\\ r\equiv r_0\mod I\endSb \Phi (r)&=
\sum \Sb r\in R\\ r\equiv wr_0\mod I\endSb \Phi (r)\\
&=\sum \Sb r\in R\\ rw^{-1}\equiv r_0\mod I\endSb \Phi (r)\\
&=\sum \Sb r'\in R\\ r'\equiv r_0\mod I\endSb \Phi (r'w)\\
&=\Phi (w)\sum \Sb r'\in R\\ r'\equiv r_0\mod I\endSb \Phi (r').
\endaligned$$
The lemma follows. \enddemo

\proclaim{Lemma 12} Suppose $R$ is a finite commutative ring with identity and
$\Phi$ is a primitive multiplicative character on $R$. Then for all 
$r_0\in R$ we have
$$\sum _{r\in R}\overline{\Phi (r)}\chi (rr_0)=
\Phi (r_0)\tau (\overline{\Phi }).$$
Also,
$|\tau (\Phi )|^2=\# R.$
\endproclaim

\demo{Proof} The first part is simply Lemma 10 when $r_0$ is a unit,
so assume otherwise (but $r_0\neq 0$) and let $I=\text{ann}(r_0)\neq\{0\}$.
Write $R/I=\{r_j+I\: j=1,\ldots ,n\}$. Then by Lemma 11
$$\aligned \sum _{r\in R}\Phi (r)\overline{\chi (rr_0)}&=\sum _{j=1}^n
\sum \Sb r\in R\\ r\equiv r_j\mod I\endSb\Phi (r)\overline{\chi (rr_0)}\\
&=\sum _{j=1}^n\overline{\chi (r_jr_0)}\sum \Sb r\in R\\ r\equiv r_j
\mod I\endSb \Phi (r)\\
&=0.\endaligned$$
Since $\Phi (r_0)=0$ by hypothesis, the first part of the lemma follows when
$r_0\neq 0$. When $r_0=0$ we have
$$\sum\Sb r\in R\endSb \overline{\Phi (r)}\chi (r0)=\overline{
\sum\Sb r\in R\endSb \Phi (r)}=0=\Phi (0)\tau (\overline{\Phi})$$
by (8).

Now since $\ker (\chi )$ doesn't contain any non-zero principal ideal
$$\aligned
\# R^{\times}|\tau (\Phi )|^2&=\sum \Sb r\in R\endSb
\left |\sum _{s\in R}\Phi (s)\chi (rs)\right |^2\\
&=\sum \Sb r\in R\endSb
\left ( \sum _{s\in R}\Phi (s)\chi (rs)\right )
\left ( \sum _{t\in R}\overline{\Phi (t)\chi (rt)}
\right )\\
&=\sum \Sb r\in R\endSb
\left ( \sum _{s\in R}\Phi (s)\chi (rs)\right )
\left ( \sum _{t\in R}\overline{\Phi (t)}\chi (-rt)\right )\\
&=\sum _{s\in R}\sum _{t\in R}
\Phi (s)\overline{\Phi (t)}
\sum \Sb r\in R\endSb\chi \big ((s-t)r\big )\\
&=\sum _{s\in R}
\Phi (s)\overline{\Phi (s)}\# R\\
&=\sum _{s\in R}|\Phi (s)|^2\# R\\
&=\# R^{\times}\# R.\endaligned$$
\enddemo

\proclaim{Lemma 13} Suppose $R$ is a finite commutative ring with identity and
$\Phi$ is a primitive character on $R$. Then for all subgroups 
$H\subseteq R$ of the additive structure of $R$, we have
$$\left |\sum _{h\in H}\Phi (h)\right |\le {\#\{u\in R^{\times}\: uH\subseteq
\ker (\chi )\}\cdot \# H\over \sqrt{\# R}}.$$
\endproclaim

\demo{Proof} By Lemma 12 
$$\Phi (h)={1\over\tau (\overline{\Phi })}\sum _{r\in R}
\overline{\Phi (r)}\chi (rh)$$
for all $h\in H$ and $|\tau (\overline{\Phi})|=\sqrt{\# R}.$ Thus
$$\aligned
\sqrt{\# R}\left |\sum _{h\in H}\Phi (h)\right |&=
\left |\sum _{h\in H}\sum _{r\in R}\overline{\Phi (r)}
\chi (rh)\right |\\
&\le \sum _{r\in R}|\overline{\Phi (r)}|\left |\sum _{h\in H}
\chi (rh)\right |\\
&=
\sum \Sb u\in R^{\times}\endSb \left |\sum _{h\in H}\chi (uh)\right |.
\endaligned$$
By (8), the inner sum here is zero unless $uH\subseteq\ker (\chi )$,
and if that is the case, then the inner sum is simply $\# H$.
\enddemo

We now return to the more specific situation where our ring is $\foc /\pia (\fc ).$
Our first goal is an estimate for what will play the role of the subgroup $H$ in
Lemma 13.

\proclaim{Lemma 14} Let $K$ be a function field and suppose
$\fc\in\Div (K)$ with $\fc >0$. Denote the canonical map $\foc\rightarrow
\foc /\pia (\fc )$ by $\theta _{\fc}$. Suppose $\fd\in\Div (K)$ is an effective
divisor with $(\fc ,\fd )=0$ and $v_0\in M(K)$ with $v_0\not\in\supp (\fc )
\cup\supp (\fd )$. Then for all integers $n$ and all proper subgroups
$G$ of the additive structure of $\foc /\pia (\fc )$ we have
$$\#\{ u\in \big ( \foc /\pia (\fc )\big )^{\times}\: u\theta _{\fc}\big (
L(nv_0-\fd )\big )\subseteq G\}\ll q_K^{\deg (\fc )+\deg (\fd )+(1-n)\deg (v_0)},
$$
where the implicit constant depends only on $K$.
\endproclaim

\demo{Proof} For notational convenience set $R=\foc /\pia (\fc )$ and set $q=q_K$.
If $\{ u\in R^{\times}\: u\theta _{\fc}\big ( L(nv_0-\fd )\big )\subseteq G\}$
is not empty, then by replacing $G$ with $u^{-1}G$ for some $u\in R^{\times},$
we see that we may assume $\theta _{\fc} \big (L(nv_0-\fd )\big )\subseteq G$.
Set
$$\aligned
n_0&=\min\{ j\in\bz\: j\deg (v_0)-\deg (\fd )\ge 2g_K-1\},\\
n_1&=\max\{ j\in\bz\: \theta _{\fc}\big (L(jv_0-\fd )\big )\subseteq G\},\endaligned$$
and $n_2\in\bz$ to be the unique solution to
$$\deg (\fc )+2g_K-1\le n_2\deg (v_0)<\deg (\fc )+2g_K-1+\deg (v_0).$$
Then $n\le n_1$ by our hypothesis above. Since $\# R=q^{\deg (\fc )}$
by Lemma 8, we see that we may assume $n\ge n_0$. Also, as in the
proof of Lemma 8 we see that both $\theta _{\fc}\big (L(n_2v_0)\big )$
and $\theta _{\fc}\big (L((n_2+n_0)v_0)\big )$ are $R$; the last equality implies that
$n_2+n_0\ge n_1$. In summary, we have $\theta _{\fc}\big (L(n_2v_0)\big )=R$ and
we may assume that $n_0\le n\le n_1\le n_2+n_0$.

Set
$$G_n=\{ r\in R\: r\theta _{\fc}\big (L(nv_0-\fd )\big )\subseteq G\}.$$
Clearly $G_n$ is a subgroup of the additive structure of $R$. We will
estimate the order of $G_n$.

Suppose $\alpha\in L\big ( (n_1-n)v_0\big )$
and $\beta\in L(nv_0-\fd )$. Then $\alpha\beta\in L(n_1v_0-\fd )$ so that
$\theta _{\fc}(\alpha\beta )\in G$ by the definition of $n_1$, whence
$\alpha\in G_n$ and $G_n\supseteq\theta _{\fc}\big (L((n_1-n)v_0)\big )$.

By definition there is an $\alpha _0
\in L\big ( (n_1+1)v_0-\fd\big )$ with $\theta _{\fc}(\alpha _0)\not\in G$.
Suppose $\alpha\in L\big ( (n_1-n_0)v_0\big )$
and write $\ordo (\alpha )=j-n_1$ for some $j\ge n_0$. By the definition of
$n_0$ we have $j\deg (v_0)-\deg (\fd)\ge 2g_K-1$. Thus, by the Riemann-Roch
Theorem
$$\aligned
\dim _{\fq}\left (L\big ( (j+1)v_0-\fd \big )/L(jv_0-\fd )\right )&=
l\big ( (j+1)v_0-\fd \big )-l(jv_0-\fd )\\
&=(j+1)\deg (v_0)-\deg (\fd )+1-g_K\\
&\qquad -\big ( j\deg (v_0)-\deg (\fd )+1-g_K\big )\\
&=\deg (v_0).\endaligned$$
Similarly,
$$\dim _{\fq}\left (L\big ( (n_1+1)v_0-\fd \big )/L(n_1v_0-\fd )\right )=
\deg (v_0).$$
Let $\beta _1,\ldots ,\beta _{\deg (v_0)}\in L\big ( (j+1)v_0-\fd \big )$
be linearly independent over $L(jv_0-\fd )$. Then we see that
 $\alpha\beta _1,\ldots ,
\alpha\beta _{\deg (v_0)}\in L\big ((n_1+1)v_0-\fd \big )$ are linearly
independent over $L(n_1v_0-\fd )$. Thus, there is a $\beta\in L\big (
(j+1)v_0-\fd \big )$ with $\alpha\beta\equiv\alpha _0\mod L(n_1v_0-\fd )$.
Since $\theta _{\fc}(\alpha _0)\not\in G$ and $\theta _{\fc}\big (L(n_1v_0-\fd )\big )
\subseteq G$ by hypotheses, we see that $\theta _{\fc}(\alpha )\not\in G_n$
whenever $j+1\le n$. In other words, if $\alpha\in 
L\big ( (n_1-n_0)v_0\big )$ and $\theta _{\fc}(\alpha )\in G_n$, then
$\alpha\in L\big ( (n_1-n)v_0\big ).$

Write $q=p^r$ where the prime $p$ is the characteristic of $K$. 
All of our groups and subgroups above
are vector spaces over $\bfp$. Since $n_2\ge n_1-n_0$, we have by the
Riemann-Roch Theorem and Clifford's Theorem
$$\aligned
\dim _{\fq}\left ( L(n_2v_0)/L\big ( (n_1-n_0)v_0\big )\right )&=
l(n_2v_0)-l\big ((n_1-n_0)v_0\big )\\
&\le (n_2-n_1+n_0)\deg (v_0)+l\big (\fw _K-(n_1-n_0)v_0\big )\\
&\le (n_2-n_1+n_0)\deg (v_0)+g_K,\endaligned$$
where $\fw _K\in\Div (K)$ is in the canonical class, so that
$$\dim _{\bfp}\left ( L(n_2v_0)/L\big ( (n_1-n_0)v_0\big )\right )\le
r(n_2-n_1+n_0)\deg (v_0)+rg_K.$$
Set $d=r(n_2-n_1+n_0)\deg (v_0)+rg_K+1$ and let $r_1,\ldots ,r_d\in G_n$.
Since $\theta _{\fc} \big (L(n_2v_0)\big )=R\supseteq G_n$, there are
$\alpha _1,\ldots ,\alpha _d\in L(n_2v_0)$ with $\theta _{\fc} (\alpha _i)=r_i$
for $i=1,\ldots ,d$. 

Now there are $a_1,\ldots ,a_d\in\bfp$, not all zero, such that
$\alpha =\sum _{i=1}^da_i\alpha _i\in L\big ( (n_1-n_0)v_0\big )$ and
$\theta _{\fc} (\alpha )\in G_n$. By what we proved above, $\alpha\in L\big (
(n_1-n)v_0\big )$. Therefore
$$\dim _{\bfp}\left ( G_n/\theta _{\fc} \big (L\big ( (n_1-n)v_0\big )\big ) \right )
\le d-1.$$
But by the Riemann-Roch Theorem and Clifford's Theorem
$$\aligned
\dim _{\fq}\left ( L\big ( (n_1-n)v_0\big )\right )&=(n_1-n)\deg (v_0)+1-g_K
+l\big (\fw _K-(n_1-n)v_0\big )\\
&\le (n_1-n)\deg (v_0)+1,\endaligned$$
so that
$$\aligned
\dim _{\bfp }(G_n)&< d+r(n_1-n)\deg (v_0)+r\\
&= r(n_2+n_0-n)\deg (v_0)+r(1+g_K)\endaligned$$
and
$$\aligned \# G_n&\ll q^{(n_2+n_0-n)\deg (v_0)}\\
&\ll q^{\deg (\fc )+\deg (\fd )+(1-n)\deg (v_0)}.\endaligned$$
\enddemo

We are now ready to state and prove our analog of the classical Polya-Vinogradov
inequality.

\proclaim{Theorem 4} Suppose $K$ is a function field 
and $\fc\in\Div (K)$ satisfies $\fc >0$. Let $\Phi$ be a
non-principal multiplicative character on the quotient ring
$\foc /\pia (\fc )$, $\fd$ be an effective divisor with $(\fc ,\fd )=0,$
and $v_0\in M(K)$ with $v_0\not\in\supp (\fc )\cup\supp (\fd )$. 
Let $\theta _{\fc} \:\foc\rightarrow\foc /\pia (\fc )$
denote the canonical map. Then for all integers $n$ with $n\deg (v_0)\ge
\deg (\fc ) +\deg (\fd ) +2g_K-1$
$$\sum\Sb r\in\theta _{\fc} ( L(nv_0-\fd ))\endSb \Phi (r)=0,$$
and for all other integers $n\ge 0$
$$\left |\sum\Sb r\in\theta _{\fc} ( L(nv_0-\fd ))\endSb \Phi (r)\right |
\ll\min\{ q_K^{n\deg (v_0)-\deg (\fd ) }, q_K^{\deg (\fc )/2+\deg (v_0)}\},$$
where the implicit constant depends only on $K$.
\endproclaim

To see the connection with the classical case, replace
the function field $K$ with the field $\bq$ and
set $v_0$ to be the place of $\bq$ corresponding to the usual Euclidean
absolute value. Then $L(nv_0)$ corresponds to the set of integers of absolute
value no greater than $n$, and the quotient ring $\foc /\pia (\fc )$ 
to $\bz /m\bz$ for some integer $m>1$. 

\demo{Proof} Set $R=\foc /\pia (\fc )$ and $q_K=q$ as above.
Whenever $n\deg (v_0)\ge\deg (\fc )+\deg (\fd )+2g_K-1$ we have
$\theta _{\fc}\big ( L(nv_0-\fd )\big )=R$ as remarked above in
the proof of Lemma 8, so that
this case is simply an application of (8). We also note that the bound
$q^{n\deg (v_0)-\deg (\fd )}$ is the trivial bound.

The ideals of $R$ are of the form
$\theta _{\fc}\big (\pia (\fa )\big )$ for effective divisors $\fa\le\fc$. 
Let $\fa$ be
minimal among such effective divisors with $\Phi (u)=\Phi (v)$ for all units 
$u,v\in R^{\times}$
with $u\equiv v\mod \theta \big (\pia (\fa )\big )$. 
(This minimal $\fa\neq 0$ since
$\Phi$ is not the principal character, and $\fa =\fc$ when $\Phi$
is primitive.) Then $\foba (\fa )\supseteq
\foc$, and if we denote by $\theta _{\fa}$ the canonical map
$\theta _{\fa}\: \foba (\fa )\rightarrow
\foba (\fa )/\pia (\fa )$, then we see that there is a primitive
character $\Phi ^*$ on $\foba (\fa )/\pia (\fa )$ with 
$$\Phi\circ \theta _{\fc}\big ( (\alpha _v)\big )=\cases
\Phi ^*\circ \theta _{\fa}\big ( (\alpha _v)\big )&\text{if
$\ord (\alpha _v)=0$ for all $v\in\supp (\fc )\setminus\supp (\fa )$,}\\
0&\text{otherwise.}\endcases$$ 

Set $\fb =\sum\Sb v\in\supp (\fc )\setminus\supp (\fa )\endSb v$. Then
$$\aligned
\left |\sum\Sb r\in\theta _{\fc} ( L(nv_0-\fd ))\endSb\Phi (r)\right |
&\le
\left |\sum\Sb \alpha\in L(nv_0-\fd )\endSb\Phi \circ\theta _{\fc}(\alpha )\right |\\
&=
\left |\sum\Sb \alpha\in L(nv_0-\fd )\\ \ord (\alpha )=0\ 
\text{all $v\in\supp (\fb )$}\endSb
\Phi ^* \circ\theta _{\fa}(\alpha )\right |\\
&=
\left |\sum\Sb \alpha\in L(nv_0-\fd )\endSb \Phi ^* \circ\theta _{\fa}(\alpha )
\sum\Sb 0\le\fe\le\fb\\ \ord (\alpha )\ge\ord (\fe)\ \text{all
$v\in\supp (\fb )$}\endSb \mu (\fe )\right |\\
&=
\left |\sum\Sb 0\le\fe\le\fb\endSb\mu (\fe )
\sum\Sb\alpha\in L(nv_0-\fd -\fe )\endSb\Phi ^*\circ\theta _{\fa}(\alpha )
\right |.\endaligned\tag 9$$

We note that there exists a character $\chi$ on $R$ with the required properties
above (see [13, Chapter IV], for instance).
Let $\fe$ be a divisor satisfying $0\le\fe\le\fb$. Then we may apply
Lemma 14 to get
$$\#\{ r\in \big (\foba (\fa ) /\pia (\fa )\big )^{\times}\: r\theta _{\fc}
\big ( L(nv_0-\fd-\fe )\big )\subseteq\ker (\chi )\}\ll q^{\deg (\fa )+
\deg (\fd +\fe )+(1-n)\deg (v_0)}.$$
Using this together with Lemma 8, Lemma 13 and the Riemann-Roch Theorem yields
$$\aligned
\left |\sum\Sb r\in\theta _{\fc}( L(nv_0-\fd-\fe))\endSb \Phi ^* (r)\right |&\le
{q^{\deg (\fa )+\deg (\fd +\fe )+(1-n)\deg (v_0)}
\cdot \#\theta _{\fc}\big ( L(nv_0 -\fd -\fe )\big )\over\sqrt{\#\foba (\fa )
/\pia (\fa )}}\\
&\ll {q^{\deg (\fa )+\deg (\fd +\fe )+(1-n)\deg (v_0)}\cdot 
q^{l(nv_0-\fd -\fe )}\over\sqrt{q^{\deg (\fa )}}}\\
&\ll {q^{\deg (\fa )+\deg (\fd +\fe )+(1-n)\deg (v_0)}\cdot 
q^{n\deg (v_0)-\deg (\fd +\fe )}\over\sqrt{q^{\deg (\fa )}}}\\
&=q^{\deg (\fa )/2 +\deg (v_0)}.\endaligned\tag 10$$
Also, as in the proof of Lemma 8
$L(nv_0-\fd-\fe )\cap\ker (\theta _{\fa})=L(nv_0-\fd-\fe -\fa )$. 
If $n\deg (v_0)
\ge\deg (\fa )+\deg (\fd +\fe )+2g_K-1$, we have $\theta _{\fa }\big (
L(nv_0-\fd-\fe )\big )=\foba (\fa )/\pia (\fa )$ and
$$\sum\Sb r\in\theta _{\fa}(L(nv_0-\fd -\fe ))\endSb \Phi ^*(r)=0$$
by (8). Hence we may assume that $\deg (nv_0-\fd -\fe -\fa )<2g_K-1$, so that
$\# L(nv_0-\fd -\fe -\fa )\ll 1$ and
$$\left |\sum\Sb \alpha\in L(nv_0-\fd -\fe )\endSb \Phi ^*\circ\theta _{\fa}
(\alpha )\right |\ll
\left |\sum\Sb r\in\theta _{\fa}(L(nv_0-\fd -\fe ))\endSb \Phi ^*(r)\right |.$$
We combine this with (9), (10) and Lemma 0 (setting $\epsilon = 1/2$ there) to get
$$\aligned
\left |\sum\Sb r\in\theta _{\fc}(L(nv_0 -\fd ))\endSb\Phi (r)\right |
&\le
\sum\Sb 0\le\fe\le\fb\endSb \left |\sum\Sb
\alpha\in L(nv_0-\fd -\fe )\endSb\Phi ^*\circ\theta _{\fa}(\alpha )\right |\\
&\ll q^{\deg (\fb )/2}q^{\deg (\fa )/2+\deg (v_0)}\\
&\le q^{\deg (\fc )/2+\deg (v_0)}.\endaligned$$
\enddemo

\head 4. The Main Estimates\endhead

In this section we state and prove our three main estimates. These are
the main tools with which we can prove results about sums involving $L$-functions.

\proclaim{Proposition 5} Suppose $K$ is a function field with $q_K=q$ odd
and let $\fc\in\Div (K)$ be an effective divisor with $\fc\not\in 2\Div (K)$.
Then for all positive integers $m$ and all $\epsilon >0$ we have
$$\left |\sum\Sb [F:K]=2\\ g_F=m,\ q_F=q\endSb \chi (F/\fc )\right |
\ll q^mq^{(\epsilon +1/4)\deg (\fc )},$$
where the implicit constant depends only on $K$ and $\epsilon$.
If $m\le\deg (\fc )/4$, then
$$\left |\sum\Sb [F:K]=2\\ g_F=m,\ q_F=q\endSb \chi (F/\fc )\right |
\ll q^{2m},$$
where the implicit constant depends only on $K$.
\endproclaim

\proclaim{Proposition 6} Suppose  $K$ is a function field with
$q_K=q$ even and
let $\fc\in\Div (K)$ be an effective divisor with $\fc\not\in 
2\Div (K)$.
Then for all positive integers $m$ and all $\epsilon >0$ we have
$$\left |\sum\Sb [F:K]=2\\ g_F=m,\ q_F=q\endSb\chi (F/\fc )\right |\ll
q^m\cdot q^{\epsilon \deg (\fc )},$$
where the implicit constant depends only on $K$ and $\epsilon$.
\endproclaim

\proclaim{Proposition 7} Let $K$ be a function field with field of 
constants $\fq$.
Fix an effective divisor $\fc\in\Div (K)$ and 
a positive integer $m$. Then for all $\epsilon >0$
$$\multline
\sum\Sb [F:K]=2\\ g_F=m,\  q_F=q\endSb\chi (F/2\fc )
=q^{2m}{2J_Kq^{3-5g_K}\over \zeta _K(2)(q-1)}
\prod _{v\in\supp (\fc)}(1+q^{-\deg (v)})^{-1}\\
+\cases
O\big (q^{(1/2+ \epsilon ) m}q^{\epsilon\deg (\fc )} )&\text{if $q$ is odd,}\\
O\big (q^{\epsilon m}q^{\epsilon\deg (\fc )} +q^m)&\text{if $q$ is
even,}\endcases\endmultline$$
where the implicit constants depend only on $\epsilon$ and $K$. 
In particular, the number $N$ of quadratic extensions $F\supset K$ with
$q_F=q$ and $g_F=m$ satisfies
$$N=q^{2m}{2J_Kq^{3-5g_K}\over \zeta _K(2)(q-1)}+
\cases O(q^{(1/2+ \epsilon )m})&\text{if $q$ is odd,}\\
O(q^m)&\text{if $q$ is even.}\endcases$$
\endproclaim

The proof of Proposition 5 requires two more intermediate results.

\proclaim{Lemma 15} Let $K$ be a function field with $q_K=q$ odd and suppose
$v_0\in M(K)$. Every quadratic extension $F$
of $K$ has a generator in $L(nv_0)$ for some positive integer $n$.
Suppose further that $\deg (v_0)$ is odd. Let $F$ be a quadratic extension of $K$
and write 
$$2n\deg (v_0)=\deg \big (\diff (F)\big )+2d$$
for some integer $d$ satisfying $0\le d<\deg (v_0)$. Let $\fa _1,\ldots ,
\fa _{J_K}$ be representatives of the divisor classes of degree $d$.
Then there is a $j$ that is uniquely determined by $v_0$ and $F$ such that $F$ 
has a generator $\omega\in L(2nv_0-2\fa _j)$; there are precisely $(q-1)/2$ such
generators $\omega$ and they all satisfy
$$\divi (\omega )=-2nv_0+2\fa _j+\diff (F).$$
\endproclaim

\demo{Proof} Suppose $F$ is a quadratic extension and let $\omega _0$
be a generator of $F$. By the Strong Approximation Theorem there is
an $\alpha\in K$ with $\ord (\alpha )=-\ord (\omega _0)$ for all
places $v\neq v_0$ with $\ord (\omega _0)<0$ and $\ord (\alpha )\ge 0$
for all other places $v\neq v_0$. Then $\omega _1=\alpha ^2\omega _0$
generates $F$ and $\ord (\omega _1)\ge 0$
for all places $v\neq v_0$ by construction; in other words, $\omega _1
\in L(nv_0)$ for some positive integer $n$. 

Now suppose $\deg (v_0)$ is odd. We have
$$\divi (\omega _1)=\ordo (\omega _1)+\diff (F)+2\fa$$
for some $\fa\ge 0$. Since $\diff (F)$ has even degree by Proposition 2,
$\ordo (\omega _1)$ must be even; write $-2n_1=\ordo (\omega _1)$.
Since $\deg (v_0)$ is odd, by the Fundamental Theorem of Arithmetic
we have $2n\deg (v_0)=\deg\big (\diff (F)\big )+2d$ for uniquely
determined integers $n>0$ and $0\le d<\deg (v_0)$. 
Then $2(n_1-n)\deg (v_0)=2\deg (\fa )-2d$, so that
$(n-n_1)v_0 +\fa$ is a divisor of degree $d$. There is an $\alpha\in K$ and a
$j$ with $\divi (\alpha )+(n-n_1)v_0+\fa =\fa _j$. Now
$F$ is generated by $\omega =\omega _1\alpha ^2$ and
$$\aligned \divi (\omega )&=\divi (\omega _1)+2\divi (\alpha )\\
&=-2n_1v_0+\diff (F)+2\fa +2\divi (\alpha )\\
&=-2nv_0+\diff (F)+2\fa _j.\endaligned$$
If $\omega '\in L(nv_0-2\fa _{j'})$ is another generator of $F$
for some $j'$, then $\omega '=\omega\beta ^2$ for some $\beta\in K$ with
$2\divi (\beta )=\divi (\beta ^2)=2\fa _{j'}-2\fa _j$, so that
$\fa _{j'}-\fa _j$ is a principal divisor and whence $j'=j$.

Finally, $\omega ,\omega '\in L(2nv_0-2\fa _j)$ generate the same
quadratic extension $F$ with $\deg \big (\diff (F)\big )=2n\deg (v_0)-2d$
if and only if $\divi (\omega ')=\divi (\omega )$ by what we have already
shown. This is the case if and only if
$\omega '/\omega =a\in\fq ^{\times}$, and $a$ must be a square since
$\omega $ and $\omega '$ generate the same extension. Since there are
$(q-1)/2$ squares in $\fq ^{\times}$, the proof is complete.
\enddemo

\proclaim{Lemma 16} Let $K$ be a function field with $q_K=q$ odd and
let $\fc\in\Div (K)$ be an effective divisor. Suppose $v_0\in M(K)$
with $v_0\not\in\supp (\fc )$. Then there are $\fa _{i,j}\in\Div (K)$
for $0\le i<\deg (v_0)$ and $1\le j\le J_K$ with the following
two properties:  i)
$\fa _{i,1},\ldots ,\fa _{i,J_K}$ are representatives of the divisor
classes of degree $i$ for all $i=1,\ldots ,\deg (v_0)-1;$  and ii)
$\ord (\fa _{i,j})\ge 0$ for all $v\neq v_0$, with equality when
$v\in\supp (\fc )$, for all $i=1,\ldots ,\deg (v_0)-1$ and
$j=1,\ldots ,J_K$.

Suppose further that $\deg (v_0)$ is odd. Then 
$$\multline {q-1\over 2}\sum\Sb [F:K]=2\\ g_F=m,\ q_F=q\endSb\chi (F/\fc )\\
=
\sum\Sb (\fb ,\fc +v_0)=0\\ \deg (\fb )\le m-2g_K+1\endSb\mu (\fb )
\sum \Sb i,j,k\\ 0\le i<\deg (v_0)\\ 1\le j\le J_K\\
k\deg (v_0)-i=m-2g_K+1-\deg (\fb )\endSb
\sum\Sb \alpha\in L(2kv_0-2\fa _{i,j})\\ \alpha\not\in L\big ( 2(k-1)v_0-2\fa _{
i,j}\big )\endSb \chic (\alpha )\endmultline$$
for all integers $m\ge 0$.
\endproclaim

\demo{Proof}
For notational convenience set $m'=m-2g_K+1$, so that
$$\sum\Sb [F:K]=2\\ g_F=m,\ q_F=q\endSb \chi (F/\fc )=\sum\Sb \fd\in\Div (K)\\
\deg (\fd )=2m'\endSb\sum\Sb [F:K]=2\\ q_F=q\\ \diff (F)=\fd\endSb
\chi (F/\fc )$$
by (2').

Fix an $i$ and let $\fa _1,\ldots ,\fa _{J_K}$ be representatives of
the divisor classes of degree $i$. By the Strong Approximation Theorem,
for each $j=1,\ldots ,J_K$ there is an $\alpha  _j\in K$ such that
$\ord (\alpha _j)=-\ord (\fa _j)$ for all places $v\neq v_0$ with $\ord (\fa _j)
<0$ or $v\in\supp (\fc )$ and $\ord (\alpha _j)\ge 0$ for all other places 
$v\neq v_0$. Then for all $j=1,\ldots ,J_K$ we have
$\ord \big (\fa _j+\divi (\alpha _j)\big )\ge 0$ for all places $v\neq v_0$, 
with equality when $v\in\supp (\fc )$. Thus for all $j=1,\ldots ,J_K$ we have
$\ord (\fa _j)\ge 0$ for
all places $v\neq v_0$, with equality if $v\in\supp (\fc )$. In this
manner, we see that there are $\fa _{i,j}$ satisfying the two properties
in the lemma.

Fix an $i$ and $j$ and suppose $k\deg (v_0)-i\ge 0$. For 
$\alpha\in L(2kv_0-2\fa _{i,j})\setminus L\big ( 2(k-1)v_0-2\fa _{i,j})$ 
we write
$$\divi (\alpha )=-2kv_0+2\fa _{i,j}+\fa '+2\fa ''$$
where $\fa '$ and $\fa ''$ are effective divisors,
$\fa '$ is square-free and $v_0\not\in\supp (\fa '')$. 
Note that $\alpha$ generates a quadratic extension
$F$ with $\diff (F)=\fa '$ (unless $\fa '=0$).
We get an idele $(\alpha _v)$ given by 
$\alpha _v=\alpha / \pi _v^{2\ord (\fa '')}$ for all $v\in M(K)$.
Then $\divi\big ( (\alpha _v)\big )=-2kv_0+2\fa _{i,j}+\fa '.$

Now suppose $\fb$ is an effective divisor with $0\le\deg (\fb )\le m'$ and
$(\fb ,\fc +v_0)=0$. Then there are uniquely determined $i$ and $k$ with
$k\deg (v_0)-i=m'-\deg (\fb )$. For any $j=1,\ldots ,J_K$ and
$\alpha\in L(2kv_0-2\fa _{i,j})\setminus L\big ( 2(k-1)v_0-\fa _{i,j}\big )$
we have
$$2\fb +\divi (\alpha )=-2kv_0+2\fa _{i,j}+\fa '+2(\fa ''+\fb )= 
\divi\big ( (\alpha _v)\big )+2\fe ,$$
where $\fe$ is an effective divisor with  $v_0\not\in\supp (\fe )$ and
$$\aligned
2\deg (\fe )&=2\deg (\fb )+\deg\big (\divi (\alpha )\big )
-\deg \big (\divi\big ( (\alpha _v)\big )\big )\\
&=2\deg (\fb )  +2k\deg (v_0)-2i-\deg (\fa ')\\
&=2m'-\deg (\fa ').\endaligned$$
For any effective divisor $\fd\in\Div (K)$, set $\iota (\fd )$ to
be the idele $(\pi _v^{\ord (\fd )})$. By construction/definition, we then have
$$\chic (\alpha )=\chic\big (\alpha\cdot\iota (2\fb )\big )=
\chic \big ( (\alpha _v)\cdot\iota (2\fe )\big ).$$

With the above in mind, using (1) yields
$$\aligned
\sum\Sb (\fb ,\fc +v_0)=0\\ \deg (\fb )\le m'\endSb \mu (\fb )&
\sum\Sb 0\le i<\deg (v_0)\\ 1\le j\le J_K\\
k\deg (v_0)-i=m'-\deg (\fb )\endSb 
\sum\Sb\alpha\in L(2kv_0-2\fa _{i,j})\\
\alpha\not\in L\big ( 2(k-1)v_0-2\fa _{i,j}\big )\endSb \chic (\alpha )\\
&=
\sum\Sb (\fb ,\fc +v_0)=0\\ \deg (\fb )\le m'\endSb\mu (\fb )
\sum\Sb 0\le i<\deg (v_0)\\ 1\le j\le J_K\\
k\deg (v_0)-i=m'-\deg (\fb )\endSb 
\sum\Sb\alpha\in L(2kv_0-2\fa _{i,j})\\
\alpha\not\in L\big ( 2(k-1)v_0-2\fa _{i,j}\big )\endSb 
\chic \big (\alpha \cdot\iota (2\fb )\big )\\
&=
\sum\Sb (\fb ,\fc +v_0)=0\\ \deg (\fb )\le m'\endSb\mu (\fb )
\sum\Sb 0\le i<\deg (v_0)\\ 1\le j\le J_K\\
k\deg (v_0)-i=m'-\deg (\fb )\endSb 
\sum\Sb\alpha\in L(2kv_0-2\fa _{i,j})\\
\alpha\not\in L\big ( 2(k-1)v_0-2\fa _{i,j}\big )\endSb \\
&\qquad\qquad\times
\sum\Sb \fb\le\fe ,\ (\fe ,v_0)=0\\ 2\deg (\fe )=2m'-\deg (\fa ')\endSb
\chic \big ( (\alpha _v)\cdot\iota (2\fe )\big )\\
&=
\sum\Sb (\fb ,\fc +v_0)=0\\ \deg (\fb )\le m'\endSb\mu (\fb )
\sum\Sb 0\le i<\deg (v_0)\\ 1\le j\le J_K\\
k\deg (v_0)-i=m'-\deg (\fb )\endSb 
\sum\Sb\alpha\in L(2kv_0-2\fa _{i,j})\\
\alpha\not\in L\big ( 2(k-1)v_0-2\fa _{i,j}\big )\endSb \\
&\qquad\qquad\times
\sum\Sb \fb\le\fe\\ 2\deg (\fe )=2m'-\deg (\fa ')\\ (\fe ,\fc +v_0)=0\endSb
\chic \big ( (\alpha _v)\big )\\
&=
\sum\Sb 0\le i<\deg (v_0)\\ 1\le j\le J_K\\
k\deg (v_0)-i\le m'\endSb 
\sum\Sb\alpha\in L(2kv_0-2\fa _{i,j})\\
\alpha\not\in L\big ( 2(k-1)v_0-2\fa _{i,j}\big )\endSb 
\chic \big ( (\alpha _v) \big )
\sum\Sb 2\deg (\fe )=2m'-\deg (\fa ')\\ (\fe ,\fc +v_0)=0\endSb
\sum\Sb 0\le\fb \le \fe \endSb \mu (\fb )\\
&=
\sum\Sb 0\le i<\deg (v_0)\\ 1\le j\le J_K\\
k\deg (v_0)-i\le m'\endSb 
\sum\Sb\alpha\in L(2kv_0-2\fa _{i,j})\\
\alpha\not\in L\big ( 2(k-1)v_0-2\fa _{i,j}\big )\\ \deg (\fa ')=2m'\endSb 
\chic \big ((\alpha _v)\big ).\endaligned$$

Now if $\alpha\in L(2kv_0-2\fa _{i,j})\setminus L\big ( 2(k-1)v_0-2\fa _{i,j}
\big )$, then $\deg (\fa ')+2\deg (\fa '')=2k\deg (v_0)-2i.$ Therefore,
if $\deg (\fa ')=2m'$ and $k\deg (v_0)-i\le m'$, we must have $k\deg (v_0)
-i=m'$ and $\fa ''=0$, so that the idele $(\alpha _v)$ is simply $\alpha$.
Thus by Proposition 4, if
$F$ denotes the quadratic extension generated by $\alpha$, we have
$\chi (F/\fc )=\chic (\alpha )$ and $\deg \big (\diff (F)\big )=2m'.$
Applying Lemma 15 finishes the proof.
\enddemo

\demo{Proof of Proposition 5} Write $m'=m-2g_K+1$ as above. 

Suppose first that $m\le \deg (\fc )/4$. The estimate in Proposition 5 in
this case is the trivial one given Proposition 7, but we can also prove it
independently by Lemma 16. Choose any place $v_0\in M(K)$ of odd degree. Then
by Lemma 16 and (0)
$$\aligned
\left |\sum\Sb [F:K]=2\\ g_F=m,\ q_F=q\endSb \chi (F/\fc )\right |&\le
\sum\Sb \fb\ge 0\\ \deg (\fb )\le m'\endSb
\sum\Sb 0\le i\le \deg (v_0)\\ 1\le j\le J_K\\ k\deg (v_0)-i=m'-\deg (\fb )
\endSb\# L(2kv_0-2\fa _{i,j})\\
&\ll
\sum\Sb \fb\ge 0\\ \deg (\fb )\le m'\endSb
\sum\Sb 0\le i\le \deg (v_0)\\ 1\le j\le J_K\\ k\deg (v_0)-i=m'-\deg (\fb )
\endSb q^{2k\deg (v_0)-2i}\\
&\ll
q^{2m'}
\sum\Sb \fb\ge 0\\ \deg (\fb )\le m'\endSb q^{-2\deg (\fb )}\\
&\ll q^{2m'}\\
&\ll q^{2m}.\endaligned$$

For the case where $m>\deg (\fc )/4$ we will again use
Lemma 16, but we must be more careful in our choice of the place $v_0$ and the divisors
$\fa _{i,j}$.

A weak corollary to the ``prime number theorem" for function fields 
(see [8, Theorem 5.12], for example) together with Lemma 0 implies that there is
a place $v_0\in M(K)$ of odd degree with $v_0\not\in\supp (\fc )$ satisfying
$q^{\deg (v_0)}\ll q^{\epsilon\deg (\fc )},$ where the implicit constant
depends only on $K$ and $\epsilon$.

Next, for
all $k$ with $k\deg (v_0)\ge 2g_K-1+\deg (\fc )$ we have
$$\aligned
\#\{\alpha\in L(kv_0)\: \ord (\alpha )=0\ \text{all $v\in\supp (\fc )$}\}&=
\sum\Sb 0\le\fb\le\fc\endSb \mu (\fb )\# L(kv_0-\fb )\\
&=q^{k\deg (v_0)+1-g_K}\sum\Sb 0\le \fb\le\fc\endSb \mu (\fb )q^{-\deg (\fb )}\\
&=q^{k\deg (v_0)+1-g_K}\prod\Sb v\in\supp (\fc )\endSb (1-q^{-\deg (v)}).
\endaligned$$
In particular, we see that for all integers $k$ with $k\deg (v_0)\ge 2g_K+
\deg (\fc )$, there is an $\alpha\in K$ with $\ordo (\alpha )=-k$ and
$\ord (\alpha )\ge 0$ for all places $v\neq v_0,$ with equality whenever
$v\in\supp (\fc )$. Hence there is a positive integer $n$ such that
we can choose divisors $\fa _{i,j}$ 
satisfying the two properties of Lemma 16 and also 
$\ordo (\fa _{i,j})= -n$ for all $0\le i<\deg (v_0)$ and $1\le j\le J_K$.
Set $\fd _{i,j}=\fa _{i,j}+nv_0$. Then $\fd _{i,j}$ is always an effective
divisor with $(\fc +v_0 ,\fd _{i,j})=0$ 
and $L(2kv_0-2\fa _{i,j})=L\big ( 2(k+n)v_0-2\fd _{i,j}\big )$.

Clearly $\chic$ induces a multiplicative character $\Phi$ on $\foc /\pia (\fc )$
and $\Phi$ is non-principal since $\fc\not\in 2\Div (K)$. Letting $\theta _{\fc}
\: \foc\rightarrow\foc /\pia (\fc )$ denote the canonical map as in the
proof of Theorem 4, we have $\ker (\theta _{\fc})\cap L\big ( 2(k+n)v_0-
2\fd _{i,j}\big )
=L\big ( 2(k+n)v_0-2\fd _{i,j}-\fc\big )$ as explained in the proof of Lemma 10.
Now either $2(k+n)\deg (v_0)-2i\ge\deg (\fc )2g_K-1$,
in which case
$$\sum\Sb \alpha\in L\big ( 2(k+n)v_0-2\fd _{i,j}\big )\endSb \chic (\alpha )=
\# L\big ( 2(k+n)v_0-2\fd _{i,j}-\fc\big )\sum\Sb r\in\foc /\pia (\fc )\endSb
\Phi (r)=0$$
by Theorem 4, or $2(k+n)\deg (v_0)-2i<\deg (\fc )+2g_K-1$,
in which case
$$\aligned
\left |
\sum\Sb \alpha\in L\big ( 2(k+n)v_0-2\fd _{i,j}\big )\endSb \chic (\alpha )
\right |&=
\# L\big ( 2(k+n)v_0-2\fd _{i,j}-\fc\big )\left |
\sum\Sb r\in\theta _{\fc}\big ( L(2(k+n)-2\fd _{i,j})\big )\endSb\Phi (r)\right |\\
&\ll\min\{ q^{2(k+n)\deg (v_0)-\deg (\fd _{i,j})},
q^{\deg (\fc )/2+\deg (v_0)}\}\\
&=\min\{ q^{2k\deg (v_0)-2i},q^{\deg (\fc )/2+\deg (v_0)}\}\\
&\ll \min\{ q^{2k\deg (v_0)-2i},q^{\deg (\fc )(\epsilon +1/2)}\}
\endaligned$$
by Theorem 4 again. We thus have
$$\aligned
\left |\sum\Sb \alpha \in L(2kv_0-2\fa _{i,j})\\ \alpha\not\in L\big (
2(k-1)v_0-2\fa _{i,j}\big )\endSb \chic (\alpha )\right |&\le
\left |\sum\Sb \alpha \in L(2kv_0-2\fa _{i,j})\endSb\chic (\alpha )\right |+
\left |\sum\Sb \alpha \in L\big ( 2(k-1)v_0-2\fa _{i,j})\endSb\chic (\alpha )
\right |\\
&\ll \min\{ q^{2k\deg (v_0)-2i},q^{\deg (\fc )(\epsilon +1/2)}\}
.\endaligned$$
Using this estimate together with Lemma 16, we get
$$\aligned
\left |\sum\Sb [F:K]=2\\ g_F=m,\ q_F=q\endSb \chi (F/\fc )\right | &\le
\left |\sum \Sb (\fb ,\fc +v_0)=0\\ m'\ge\deg (\fb )\ge m'-\deg (\fc )/4\endSb
\mu (\fb )
\sum \Sb 0\le i<\deg (v_0)\\ 1\le j\le J_K\\ k\deg (v_0)-i=m'-\deg (\fb )\endSb
\sum\Sb \alpha\in L(2kv_0-2\fa _{i,j})\\ \alpha\not\in L\big ( 2(k-1)v_0-2
\fa _{i,j}\big )\endSb \chic (\alpha )\right |\\
&+
\left |\sum \Sb (\fb ,\fc +v_0)=0\\ 0\le\deg (\fb )< m'-\deg (\fc )/4\endSb
\mu (\fb )
\sum \Sb 0\le i<\deg (v_0)\\ 1\le j\le J_K\\ k\deg (v_0)-i=m'-\deg (\fb )\endSb
\sum\Sb \alpha\in L(2kv_0-2\fa _{i,j})\\ \alpha\not\in L\big ( 2(k-1)v_0-2
\fa _{i,j}\big )\endSb \chic (\alpha )\right |\\
&\le
\sum \Sb (\fb ,\fc +v_0)=0\\ m'\ge\deg (\fb )\ge m'-\deg (\fc )/4\endSb
\sum \Sb 0\le i<\deg (v_0)\\ 1\le j\le J_K\\ k\deg (v_0)-i=m'-\deg (\fb )\endSb
\left |
\sum\Sb \alpha\in L(2kv_0-2\fa _{i,j})\\ \alpha\not\in L\big ( 2(k-1)v_0-2
\fa _{i,j}\big )\endSb \chic (\alpha )\right |\\
&\qquad +
\sum \Sb (\fb ,\fc +v_0)=0\\ 0\le\deg (\fb )< m'-\deg (\fc )/4\endSb
\sum \Sb 0\le i<\deg (v_0)\\ 1\le j\le J_K\\ k\deg (v_0)-i=m'-\deg (\fb )\endSb
\left |
\sum\Sb \alpha\in L(2kv_0-2\fa _{i,j})\\ \alpha\not\in L\big ( 2(k-1)v_0-2
\fa _{i,j}\big )\endSb \chic (\alpha )\right |\\
&\ll 
\sum \Sb (\fb ,\fc +v_0)=0\\ m'\ge\deg (\fb )\ge m'-\deg (\fc )/4\endSb
q^{2m'-2\deg (\fb )}\\
&\qquad +
\sum \Sb (\fb ,\fc +v_0)=0\\ 0\le\deg (\fb )< m'-\deg (\fc )/4\endSb
q^{\deg (\fc )(\epsilon +1/2)}\\
&\ll q^{2m'}q^{-m'+\deg (\fc )/4}+q^{m'-\deg (\fc )/4}q^{\deg (\fc )(
\epsilon +1/2)}\\
&\ll q^{m+(\epsilon +1/4)\deg (\fc )}.\endaligned$$
\enddemo

We now turn to the case of characteristic two and Proposition 6. We first note
that by Lemma 9 it suffices to prove Proposition 6 under the assumption that
$\deg (\fc )$ is even. Before we can finish the proof, we need
several intermediate results.

\proclaim{Lemma 17} Suppose $K$ is a function field with $q_K=q$ even and
$\fc$ is an effective divisor with 
$\fc\not\in 2\Div (K)$. Then the function $\chic$ is a non-trivial character
on $\foc$. Further, $L(0)\subseteq\ker (\chic )$ if and only if
$\deg (\fc )$ is even. For each  place $v\not\in\supp (\fc )$, the supremum
$\sup \{ n\: L(nv)\subseteq\ker (\chic )\}$
is an achieved maximum; call it $n_v$ and set $n_v=0$ for all 
$v\in\supp (\fc )$.
For any effective divisor $\fa$ with $(\fa ,\fc )=0$, $L(\fa )\subseteq
\ker (\chic )$ only if $\ord (\fa )\le n_v$ for all $v\in M(K)$.
Finally, if $\deg (\fc )$ is even then we get an effective divisor
$$\fcchi =\sum\Sb v\in M(K)\endSb n_vv,$$
and for all effective divisors $\fa$ with $(\fa ,\fc )=0,$
$L(\fa )\subseteq\ker (\chic )$ only if $\fa\le\fcchi$.
\endproclaim

\demo{Proof} Clearly $\chic$ is a character on $\foc$ from the definitions,
with $\ker (\chic )\supseteq\pia (\fc )$.
Since $\chi _v$ is certainly non-trivial for all places $v\in M(K)$, we see that
$\chic$ is non-trivial whenever $\fc\not\in 2\Div (K)$ by the definitions.

Next, $L(0)=\fq$. Clearly all $a\in\fq$ of the form $a=b^2+b$ for
some $b\in\fq$ satisfy $\chi _v(a)=1$ for all places $v\in M(K)$.
Thus by Lemma 11 we see that $L(0)=\fq\subseteq\ker (\chi _{\fc})$
if and only if $\deg (\fc )$ is even.

Note that $L(-v)=\{ 0\}$ for all places $v\in M(K),$ so 
that $n_v\ge -1$ always. If $v\not\in\supp (\fc )$, then
$L(nv)\not\subseteq\ker (\chic )$ for all positive integers $n$ with
$n\deg (v)\ge\deg (\fc _1)+2g_K-1$ by Lemma 10, since $\chic$ is non-trivial
with $\ker (\chic )\supseteq\pia (\fc )$. Thus $n_v$ is always an achieved
maximum.

Now suppose $\fa$ is an effective divisor with $(\fa ,\fc )=0$ and
$L(\fa )\subseteq\ker (\chic )$. Let $v\in\supp (\fa )$. Then since
$\fa$ is an effective divisor
$L(\ord (\fa )v)\subseteq L(\fa )\subseteq\ker (\chic )$, so that
$\ord (\fa )\le n_v$. 

Finally, $n_v\le 0$ for all but finitely many places $v\in M(K)$ (indeed, for
all places of degree at least $\deg (\fc _1)+2g_K-1$), whence $n_v=0$ for
all but finitely many places whenever $L(0)\subseteq\ker (\chic )$.
\enddemo

\proclaim{Lemma 18} Suppose $K$ is a function field with $q_K=q$ even 
and $\fd$ is an
effective divisor with $\deg (\fd _2)\ge 2g_K-1$.
Then every $F\in S(\fd )$ has ${q^{\deg (\fd _2)}\over 2}$
generators in $L'(\fd _1+2\fd _2)$ and $N(\fd )=2\phi (\fd )$.
In particular, for all effective divisors $\fc$ with $(\fc ,\fd )=0$ we have
$$\aligned {q^{\deg (\fd _2)}\over 2}
\sum\Sb [F:K]=2,\ q_F=q\\ \diff (F)=2\fd\endSb \chi (F /\fc )&=
\sum\Sb \omega\in L'(\fd _1+2\fd _2)\endSb
\chic (\alpha )\\
&=\sum\Sb 0\le\fa\le\fd _1+2\fd _2\endSb \mu (\fa )
\sum\Sb \alpha\in L(\fd _1+2\fd _2-\fa )\endSb \chic (\alpha ).\endaligned$$
\endproclaim

\demo{Proof} The first part follows directly from Lemma 5 and Proposition 3.
The last part follows from the first part, Proposition 4
and a simple application of M\"obius inversion using (1).
\enddemo

\proclaim{Lemma 19} Suppose $K$ is a function field with $q_K=q$ even 
and $\fc$ is an effective divisor with
$\fc\not\in 2\Div (K)$. 
Suppose $\fd$ is an effective divisor with $(\fc ,\fd )=0$ and
$\deg (\fd _2)\ge 2g_K-1.$ If $\deg (\fc )$ is even and
$2\fd _2\not\le\fcchi$, then 
$$\sum\Sb [F:K]=2,\ q_F=q\\ \diff (F)=2\fd\endSb \chi (F/\fc )=0.$$
\endproclaim

\demo{Proof} Since $\deg (\fd _2)\ge 2g_K-1$ we may use Lemma 18. 
Assume $\deg (\fc )$ is even and $2\fd _2\not\le\fcchi$. We 
note that if $0\le\fa\le\fd _1+2\fd _2$, then $\mu (\fa )=0$
unless $\fa\le\fd _1$. Moreover, if $\fa\le\fd _1,$ then
$\fd _1+2\fd _2-\fa \ge 2\fd _2$ and by Lemma 17 we have 
$L(2\fd _2)\subseteq L(\fd _1+2\fd _2
-\fa )\not\subseteq\ker (\chic )$ since $2\fd _2\not\le\fcchi$. We thus have
$$\aligned
{q^{\deg (\fd _2)}\over 2}\sum\Sb [F:K]=2,\ q_F=q_K\\
 \diff (F)=2\fd \endSb\chi (F/\fc )&=
\sum\Sb 0\le\fa\le\fd _1+2\fd _2\endSb\mu (\fa )\sum\Sb \alpha\in
L(\fd _1+2\fd _2-\fa )\endSb \chic (\alpha )\\
&=
\sum\Sb 0\le\fa\le\fd _1\endSb\mu (\fa )\sum\Sb \alpha\in
L(\fd _1+2\fd _2-\fa )\endSb \chic (\alpha )\\
&=
\sum\Sb 0\le\fa\le\fd _1\endSb\mu (\fa )\cdot 0\\
&=0.\endaligned$$
\enddemo

\proclaim{Lemma 20} Suppose $K$ is a function field with $q_K=q$ even and
let $\fc,\ \fd$ be effective divisors with
$(\fc ,\fd )=0$ and $\fc\not\in 2\Div (K)$. 
Let $v_0\in M(K)$ with $v_0\not\in \supp (\fc )\cup\supp (\fd )$ and 
$\deg (v_0)\ge 2g_K-1$. If $\fd >0$, then
$$\sum _{i=0}^{[n_{v_0}/2]+1}
\sum\Sb [F:K]=2,\ q_F=q\\ \diff (F)=2\fd+2iv_0\endSb \chi (F/\fc )=
0.$$
If $\fd =0$, then
$$\sum _{i=1}^{[n_{v_0}/2]+1}
\sum\Sb [F:K]=2,\ q_F=q\\ \diff (F)=2iv_0\endSb \chi (F/\fc )=
-1-(-1)^{\deg (\fc )}.$$
\endproclaim

\demo{Proof} Suppose $\fd >0$. As in the proof of Proposition 3, we may
choose a positive integer $n$ such that all $F\in S(\fd )\cup S(\fd +v_0)$
have a generator in
$$S=\{ \omega\in L(\fd _1+2\fd _2+2nv_0\: \ord (\omega )=1-2\ord (\fd )\ \text{for
all $v\in\supp (\fd )$}\}.$$
We may assume without loss of generality that $n\ge [n_{v_0}]/2+1$.
From the proof of Proposition 3 we have
$$ {2\over q^{l(\fd _2+nv_0)}}
\sum\Sb\omega\in S\endSb \chic (\omega )=
\sum _{i=0}^n\sum\Sb [F:K]=2,\ q_F=q\\
\diff (F)=2\fd +2iv_0\endSb\chi (F/\fc ).\tag 11$$
On the other hand, 
by M\"obius inversion exactly as in the proof of Lemma 4 we have
$$\sum\Sb\omega\in S\endSb \chic (\omega )=\sum\Sb 0\le\fa
\le\fd _1+2\fd _2\endSb \mu (\fa )\sum\Sb\alpha\in L(\fd _1+2\fd _2
+2nv_0-\fa )\endSb \chic (\alpha ).$$
Since $L(\fd _1+2\fd _2+2nv_0-\fa)\supseteq L(2nv_0)$ for all
divisors $\fa$ with $0\le\fa\le\fd _1+2\fd _2$ and $2n>n_{v_0}$,
we have $L(\fd -1+2\fd _2+2nv_0-\fa )\not\subseteq\ker (\chic )$, so that
$$\sum\Sb\omega\in S\endSb \chic (\omega )=0\tag 12$$
by (8).

Now suppose $n\ge i>[n_{v_0}/2]+1$. We note that $(\fd +iv_0)_2=
\fd _2+(i-1)v_0$. Since $\deg \big ( (i-1)v_0\big )\ge\deg (v_0)\ge 2g_K-1$
and $2(i-1)>n_{v_0}$, we get
$$\sum\Sb [F:K]=2,\ q_F=q_K\\ \diff (F)=2\fd +2iv_0\endSb\chi (F/\fc )=0
\tag 13$$
by Lemma 19. The case where $\fd >0$ now follows from (11)-(13).

Suppose now that $\fd =0$. We still have (12) and (13), but now
$S=L(2nv_0)$ and (as remarked in the proof of Proposition 3) contains
$q^{l(nv_0)}/2$ elements $\alpha$ of the form $\alpha =\beta ^2+\beta$
and $q^{l(nv_0)}/2$ generators of the field $K(a)$, where $a\in\fq$
is as in Lemma 11. We obviously have $\chic (\beta ^2+\beta )=1$, and
$\chic (a)=(-1)^{\deg (\fc )}$ by Lemma 11. Thus we replace (11) with
$$ {2\over q^{l(nv_0)}}
\sum\Sb\alpha\in S\endSb \chic (\alpha )=1+(-1)^{\deg (\fc )}+
\sum _{i=1}^n\sum\Sb [F:K]=2,\ q_F=q\\
\diff (F)=2\fd +2iv_0\endSb\chi (F/\fc ).$$
The lemma follows.
\enddemo

\proclaim{Lemma 21} Suppose $K$ is a function field with $q_K=q$ even 
and $\fc$ is an effective
divisor with $\fc\not\in 2\Div (K)$. Let $\fd$ be an effective divisor
with $(\fc ,\fd )=0$ of the form 
$$\fd =\fd '+v_1+\cdots +v_m,$$
where all places $v\in\supp (\fd ')$ satisfy $\deg (v)<2g_K-1$
and the $v_1,\ldots ,v_m$ are distinct places with $\deg (v_i)\ge 2g_K-1$
for all $i=1,\ldots ,m$. If $\fd '>0$, then
$$\aligned
\sum\Sb [F:K]=2,\  q_F=q\\ \diff (F)=2\fd\endSb\chi (F/\fc )&=
(-1)^m\sum\Sb [F:K]=2,\ q_F=q\\ \diff (F)=2\fd '\endSb \chi (F/\fc )\\
&+\sum _{i=1}^m (-1)^{m-i-1}\sum_{j=2}^{[n_{v_i}/2]+1}\sum\Sb [F:K]=2,\
q_F=q\\ \diff (F)=2\fd '+2v_1+\cdots +2v_{i-1}+2jv_i\endSb\chi (F/\fc ).
\endaligned$$
In particular, if $n_{v_m}\le 0$, then
$$\sum\Sb [F:K]=2,\ q_F=q\\ \diff (F)=2\fd \endSb \chi (F/\fc )=
-\sum\Sb [F:K]=2,\ q_F=q\\ \diff (F)=2(\fd -v_m) \endSb \chi (F/\fc ).$$

If $\fd '=0$, then
$$\aligned
\sum\Sb [F:K]=2,\  q_F=q\\ \diff (F)=2\fd\endSb\chi (F/\fc )&=
(-1)^m\big (1+(-1)^{\deg (\fc )}\big )\\
&+\sum _{i=1}^m (-1)^{m-i-1}\sum_{j=2}^{[n_{v_i}/2]+1}\sum\Sb [F:K]=2,\
q_F=q\\ \diff (F)=2v_1+\cdots +2v_{i-1}+2jv_i\endSb\chi (F/\fc ).
\endaligned$$
In particular, if $n_{v_m}\le 0$, then
$$\sum\Sb [F:K]=2,\ q_F=q\\ \diff (F)=2\fd \endSb \chi (F/\fc )=
-\sum\Sb [F:K]=2,\ q_F=q\\ \diff (F)=2(\fd -v_m) \endSb \chi (F/\fc )$$
if $m>1$, and if $m=1$
$$\sum\Sb [F:K]=2,\ q_F=q\\ \diff (F)=2v_1 \endSb \chi (F/\fc )=
-1-(-1)^{\deg (\fc )}.$$
\endproclaim

\demo{Proof} First assume $\fd '>0$.
We prove the lemma by induction on $m$. The case $m=1$ is
just Lemma 20 with the sums reorganized. Now assume $m>1$ and the
lemma is true for $m-1$. By Lemma 20 we have 
$$\aligned
\sum\Sb [F:K]=2,\ q_F=q\\ \diff (F)=2\fd\endSb\chi (F/\fc )&=
-\sum\Sb [F:K]=2,\ q_F=q\\ \diff (F)=2\fd '+2v_1+\cdots +2v_{m-1}\endSb
\chi (F/\fc )\\
&-\sum _{j=2}^{[n_{v_m}/2]+1}
\sum\Sb [F:K]=2,\ q_F=q\\ \diff (F)=2\fd '+2v_1+\cdots +2v_{m-1}
+2jv_m\endSb\chi (F/\fc ).\endaligned$$
Now by the case $m-1$,
$$\multline
\sum\Sb [F:K]=2,\ q_F=q\\ \diff (F)=2\fd '+2v_1+\cdots +2v_{m-1}\endSb
\chi (F/\fc )=
(-1)^{m-1}\sum\Sb [F:K]=2,\ q_F=q\\ \diff (F)=2\fd '\endSb\chi (F/\fc )\\
-\sum _{i=1}^{m-1}(-1)^{m-i-1}\sum _{j=0}^{[n_{v_i}/2]+1}
\sum\Sb [F:K]=2,\ q_F=q\\ \diff (F)=2\fd '+2v_1+\cdots +2v_{i-1}+2jv_i
\endSb\chi (F/\fc ).\endmultline$$
This proves the lemma when $\fd '>0$. The proof when $\fd '=0$ is similar.
\enddemo

\proclaim{Lemma 22} Suppose $K$ is a function field with $q_K=q$ even 
and $\fc$ is an effective divisor
with $\fc\not\in 2\Div (K)$. Suppose $\deg (\fc )$ is even. Then
$$\deg (\fcchi )\le \deg (\fc _1)+g_K\#\supp (\fcchi ).$$ 
\endproclaim

\demo{Proof} 
First, since $\chic$ is a character and $L(n_vv)\subseteq\ker (\chic )$
for all places $v\in M(K)$ by definition, we have
$$\sum\Sb v\in\supp (\chic )\endSb L(n_vv)\subseteq\ker (\chic ).$$
Next, $\chic$ is non-trivial by Lemma 17 and $\pia (\fc )\subseteq\ker (\chic )$
as remarked in the proof of Lemma 17. Thus, the image of $\ker (\chic )$
under the canonical map $\foc\mapsto\foc /\pia (\fc )$  is a proper subspace.
Now by Lemma 10 we have
$$\deg (\fc _1)=\dim _{\fq}\foc /\pia (\fc )>\dim _{\fq }
\left (\sum\Sb v\in\supp (\chic )\endSb L(n_vv)\right ).\tag 14$$
We claim that for any positive integer $m$ and effective divisors
$\fa _1,\ldots ,\fa _m$ with $(\fa _i,\fa _j)=0$ for all $i\neq j$,
we have
$$\dim _{\fq}\left (\sum_{i=1}^m L(\fa _i)\right )
\ge\deg\left  (\sum _{i=1}^m
\fa _i\right )+1-mg_K.\tag 15$$
Note that the lemma
follows from (14) and (15).
We prove (15) by induction on $m$.

The case $m=1$ in (15) follows immediately from the Riemann-Roch Theorem.
Now suppose $m>1$. Note that 
$$\left ( \sum _{i=1}^{m-1}L(\fa _i)\right )\cap  L(\fa _m)=L(0)=\fq.$$
Now by the case $m-1$ of (15) and the Riemann-Roch Theorem we get
$$\aligned
\dim _{\fq}\left (\sum _{i=1}^mL(\fa _i)\right )&=
\dim _{\fq}\left (\sum _{i=1}^{m-1}L(\fa _i)\right )+l(\fa _m)-l(0)\\
&\ge\deg \left (\sum_{i=1}^{m-1}\fa _i\right )+1-(m-1)g_K+\deg (\fa _m)
+1-g_K-1\\
&=\deg\left (\sum _{i=1}^m\fa _i\right )+1-mg_K.\endaligned$$
\enddemo

\demo{Proof of Proposition 6} By Proposition 1,
$$\sum\Sb [F:K]=2\\ g_F=m,\ q_F=q\endSb\chi (F/\fc )=
\sum\Sb \fd >0\\ \deg (\fd )=m+1-g_K\endSb\sum\Sb [F:K]=2,\ q_F=q\\
\diff (F)=2\fd\endSb\chi (F/\fc ).\tag 16$$	

Suppose $\deg (\fc )$ is even.  For all effective divisors
$\fd$ with $(\fd ,\fc )=0$, $\deg (\fd _2)\ge 2g_K-1$ and $2\fd _2\not
\le\fcchi$, we have
$$\sum\Sb [F:K]=2,\ q_F=q\\ \diff (F)=2\fd\endSb\chi (F/\fc )=0\tag 17$$
by Lemma 19. So suppose either $\deg (\fd _2)<2g_K-1$ or $2\fd _2\le\fcchi$.
This time we write $\fd =\fa +\fb$ where
$$\fa =\sum\Sb v\in\supp (\fd )\cap\supp (\fcchi )\endSb \ord (\fd )v,
\qquad
\fb =\sum\Sb v\in\supp (\fd )\\ v\not\in\supp (\fcchi )\endSb 
\ord (\fd )v.$$
Write $\fb =\fb '+v_1+\cdots +v_n$ as in Lemma 21. We note that the places
$v_i$, $1\le i\le n$ are distinct since otherwise we would have
$\deg (\fd _2)\ge 2g_K-1$ and $\fd _2\not\le\fcchi$. Also $n_{v_i}=0$
for all $i=1,\ldots ,n$ since $(\fb ,\fcchi )=0$ by construction. 
Since $\deg (\fc )$ is even, Lemma 21 implies that
$$\sum\Sb [F:K]=2,\ q_F=q\\ \diff (F)=2\fd\endSb\chi (F/\fc )=\cases
\pm\sum\Sb [F:K]=2,\ q_F=q\\ \diff (F)=2\fa +2\fb '\endSb\chi (F/\fc )&
\text{if $\fa +\fb '\neq 0$,}\\
-2&\text{if $\fa +\fb '=0$.}\endcases\tag 18$$
We next note that for all places $v\in\supp (\fb ')$ we have
$\ord (\fb ')<2g_K$, since otherwise we would have $\deg (\fd _2)\ge 2g_K-1$
and $\fd _2\not\le\fcchi$. Thus, if we set
$$\fe =\sum\Sb v\in M(K)\\ \deg (v)<2g_K-1\endSb 2g_Kv,$$
then $\fa +\fb '\le\fcchi +\fe$.

Now  using Lemma 0 together with Lemma 22 and $\deg (\fe )\ll 1$,
we see that for all $\epsilon >0$
the number $N$ of effective divisors $\fa +\fb '\le\fcchi+\fe$ satisfies
$N\ll q^{\epsilon\deg (\fc _1)}\le 
q^{\epsilon\deg (\fc )}$. Here
the implicit constant depends only on $\epsilon$ and $K$.
Also, we clearly have
$$\left |
\sum\Sb [F:K]=2,\ q_F=q\\ \diff (F)=2\fa +2\fb '\endSb\chi (F/\fc )\right |\le
\sum\Sb [F:K]=2,\ q_F=q\\ \diff (F)=2\fa +2\fb '\endSb 1
< q^{\deg (\fa +\fb ')}\tag 19$$
by Proposition 3. Now by (16)-(19)
$$\aligned \left |\sum\Sb [F:K]=2\\ q_F=q,\ g_F=m\endSb \chi (F/\fc )\right |
&\le\sum _{i=0}^{m+1-g_K}\sum \Sb \fb \ge 0\\ \deg (\fb )=i\endSb
\left |\sum \Sb [F:K]=2,\ q_F=q\\ \diff (F)=2\fa +2\fb '\\ \deg (\fa +\fb ')=
m+1-g_K-i\endSb\chi (F/\fc )
\right |\\
&\ll\sum _{i=0}^{m+1-g_K}\sum \Sb \fb \ge 0\\ \deg (\fb )=i\endSb
q^{\epsilon\deg (\fc )}q^{m+1-g_K-i}\\
&\ll q^{m}q^{\epsilon \deg (\fc )}.\endaligned$$
\enddemo

\demo{Proof of Proposition 7} We first note that the last part of 
Proposition 7 follows from the
first part by setting $\fc =0$. We also note that
by Proposition 1, $\deg \big (\diff (F)\big )=2m-4g_K+2$ if $g_F=m$,
 $\chi (F/2\fc )=1$ if $\big (\diff (F),\fc \big )=0$ and 
$\chi (F/2\fc )=0$ if $\big (\diff (F),\fc \big )\neq 0$. 
It will prove convenient to temporarily set
$$\chi _{\fc}(\fd )=\cases 1&\text{if $(\fc ,\fd )=0$,}\\
0&\text{otherwise.}\endcases$$

For the case where $q$ is odd 
set $m'=m-2g_K+1$ and choose a place $v_0\not\in\supp (\fc )$ of odd degree 
and divisors $\fa _{i,j}$ as in the proof of Proposition 5 above.

For a fixed $i,\ j$ and $k$,
$$\sum\Sb\alpha\in L(2kv_0-2\fa _{i,j})\endSb\chi _{2\fc} (\alpha )=
\#\{\alpha\in L(2kv_0-2\fa _{i,j})\: \ord (\alpha )=0\ \text{all
$v\in\supp (\fc )$}\}.$$
For any divisor $\fa\in\Div (K)$ with $\supp (\fa )\cap\supp (\fc )=\emptyset$,
one readily verifies that
$$\#\{\alpha\in L(\fa )\: \ord (\alpha )=0\ \text{all $v\in\supp (\fc )$}\}
=\sum\Sb 0\le\fd\le\fc\endSb\mu (\fd )\# L(\fa -\fd ).$$
If $\deg (\fa )\ge 2g_K-1+\deg (\fc )$, then
$$\aligned\sum\Sb 0\le \fd\le\fc\endSb \mu (\fd )\# L(\fa -\fd )&=
\sum\Sb 0\le\fd\le\fc\endSb\mu (\fd )q^{\deg (\fa -\fd ) +1-g_K}\\
&q^{\deg (\fa )+1-g_K}\sum\Sb 0\le\fd\le\fc\endSb\mu (\fd )q^{-\deg (\fd )}\\
&=
q^{\deg (\fa )+1-g_K}\prod\Sb v\in\supp (\fc )\endSb (1-q^{-\deg (v)}).
\endaligned$$
If $\deg (\fa )<2g_K-1+\deg (\fc )$, then
$$\aligned
\left |
\sum\Sb 0\le\fd\le\fc\endSb\mu (\fd )\big (\# L(\fa -\fd )-q^{\deg (
\fa -\fd )+1-g_K}\big )\right |&\ll
\sum\Sb 0\le\fd\le\fc\endSb 1\\
&\ll q^{\epsilon\deg (\fc )}\endaligned$$
by Lemma 0 and the Riemann-Roch Theorem. 
Thus, whenever $2k\deg (v_0)-2i\ge \deg (\fc )+2\deg (v_0)+2g_K-1$ we have
$$\sum\Sb \alpha\in L(2kv_0-2\fa _{i,j})\\ \alpha\not\in L\big (
2(k-1)v_0-\fa _{i,j}\big )\endSb\chi _{2\fc} (\alpha )=
q^{2k\deg (v_0)-2i+1-g_K}(1-q^{-2\deg (v_0)})\prod\Sb v\in\supp (\fc )\endSb
(1-q^{-\deg (v)}),$$
and if $2k\deg (v_0)-2i<\deg (\fc )+2\deg (v_0)+2g_K-1$
$$\multline \left |
\sum\Sb \alpha\in L(2kv_0-2\fa _{i,j})\\ \alpha\not\in L\big (
2(k-1)v_0-\fa _{i,j}\big )\endSb\chi _{2\fc} (\alpha )-
q^{2k\deg (v_0)-2i+1-g_K}(1-q^{-2\deg (v_0)})\prod\Sb v\in\supp (\fc )\endSb
(1-q^{-\deg (v)})\right |\\
\ll q^{\epsilon /2\deg (\fc )}.\endmultline$$

For notational convenience, temporarily set $c=\deg (\fc )+2\deg (v_0)+2g_K-1$.
Combining the above together with Lemma 16 yields 
$$\multline
{q-1\over 2}\sum\Sb [F:K]=2\\ g_F=m,\ q_F=q\endSb  \chi (F/2\fc )\\
=J_K\sum\Sb (\fb ,\fc +v_0)=0\\ \deg (\fb )\le m'\endSb\mu (\fb )
q^{2m'-2\deg (\fb )+1-g_K}(1-q^{-2\deg (v_0)})\prod\Sb v\in\supp (\fc )\endSb
(1-q^{-\deg (v)})\\
+O\left (
\sum\Sb m'-c/2<n\le m'\endSb q^{\epsilon /2\deg (\fc)}
\left |\sum\Sb (\fb ,\fc +v_0)=0\\ \deg (\fb )=n\endSb\mu (\fb )\right |\right ).
\endmultline\tag 20$$
Now via the Euler product and the definition of $L_K$,
$$\aligned
\sum\Sb \fb \ge 0\\ (\fb ,\fc +v_0)=0\endSb\mu (\fb )q^{-s\deg (\fb )}&=
\sum\Sb\fb\ge 0\endSb \mu (\fb )\chi _{\fc +v_0}(\fb )q^{-s\deg (\fb )}\\
&={\prod \Sb v\in M(K)\endSb (1-q^{-s\deg (v)})\over
\prod \Sb v\in\supp (\fc +v_0)\endSb (1-q^{-s\deg (v)})}\\
&={(1-q^{-s})(1-q^{1-s})\over L_K(q^{-s})}
\prod\Sb v\in\supp (\fc +v_0)\endSb (1-q^{-s\deg (v)})^{-1}
.\endaligned\tag 21$$
We note by Lemma 0 that
$$\left |\prod\Sb v\in \supp (\fc +v_0)\endSb (1-q^{-s\deg (v)})^{-1}\right |
\ll \left |\prod\Sb v\in \supp (\fc )\endSb (1-q^{-s\deg (v)})^{-1}\right |
\ll q^{\epsilon /2\deg (\fc )}\tag 22$$
for all $s\in\bc$ with $\Re (s)\ge 1/2$. 

Combining (21) and (22) yields
$$\left |\sum \Sb (\fb ,\fc +v_0)=0\\ \deg (\fb )=n\endSb \mu (\fb )
\right |\ll q^{(1/2+\epsilon )n}q^{\epsilon /2 \deg (\fc )}\tag 23$$
for all $n\ge 0$. Also by (21) (with $s=2$) and (22), we have
$$\aligned 
\sum\Sb (\fb ,\fc +v_0)=0\\ \deg (\fb )\le m'\endSb \mu (\fb )q^{-2
\deg (\fb )}&=
\zetak (2)^{-1}(1-q^{-2\deg (v_0)})^{-1}\\
&\qquad\times \prod _{v\in\supp (\fc )}
(1-q^{-\deg (v)})^{-1}(1+q^{-\deg (v)})^{-1}\\
&\qquad\qquad +O(q^{(-3/2+\epsilon )m'}q^{\epsilon\deg (\fc )}).\endaligned\tag 24$$
The case where $q$ is odd follows from (20), (23) and (24).

We now turn to the case of even characteristic.  
For $s\in\bc$ with $\Re (s)>1$ we set
$$\aligned f(s)&=\sum\Sb \fd\ge 0\endSb q^{-(s+1)\deg (\fd )}\phi (\fd )
\chi _{\fc }(\fd )\\
&=\sum\Sb \fd\ge 0\endSb q^{-s\deg (\fd )}q^{-\deg (\fd )}\phi (\fd )
\chi _{\fc }(\fd ).\endaligned$$
Since
$q^{-\deg (\fa +\fb )}\phi (\fa +\fb )\chi _{\fc }(\fa +\fb)
=q^{-\deg (\fa )}\phi (\fa )\chi _{\fc }(\fa )q^{-\deg (\fb )}
\phi (\fb )\chi _{\fc }(\fb )$ whenever $(\fa ,\fb )=0,$
we have
$$\aligned 
f(s)&=
\prod _{v\in M(K)}\big ( 1+q^{-\deg (v)}\phi (v)\chi _{\fc }(v)q^{-s\deg (v)}+
q^{-2\deg (v)}\phi (2v)\chi _{\fc }(2v)q^{-2s\deg (v)}+\cdots \big )\\
&=
{\prod _{v\in M(K)}\big ( 1+(1-q^{-\deg (v)})q^{-s\deg (v)}+(q^{\deg (v)}-
1)q^{-2s\deg (v)}+\cdots \big )\over
\prod _{v\in\supp (\fc )}\big ( 1+(1-q^{-\deg (v)})q^{-s\deg (v)}+(q^{\deg (v)}-
1)q^{-2s\deg (v)}+\cdots \big )}\\
&=
{\prod _{v\in M(K)}\big ( 1-q^{-s\deg (v)}\big )^{-1}\big ( 1-q^{-(s+1)\deg (v)}
\big )\over
\prod _{v\in\supp (\fc )}\big ( 1-q^{-s\deg (v)}\big )^{-1}
\big ( 1-q^{-(s+1)\deg (v)}\big )}\\
&=
{\zeta _K(s)\over\zeta _K(s+1)}\prod _{v\in\supp (\fc )}
\big ( 1-q^{-s\deg (v)}\big )
\big ( 1-q^{-(s+1)\deg (v)}\big )^{-1}.\endaligned$$
We see that $f(s)$ is holomorphic on any set $\{ s\in\bc\:
\Re (s)\ge \epsilon -1/2, \ -\pi /\log q\le\Im (s)<\pi /\log q \}$ 
except for a simple pole at $s=1$. 

Using a standard Tauberian argument (see [8, Theorem 17.1], for example)
we have
$$\sum\Sb\fd\ge 0\\ \deg (\fd )=j\endSb q^{-j}\phi (\fd )\chi _{\fc }
(\fd )=r\log (q) q^j+O(Mq^{(\epsilon -1)j}),\tag 25$$
where $r$ is the residue of $f(s)$ at $s=1$,
$M$ is the maximum of $|f(s)|$ on the set of $s$ with $\Re (s)=\epsilon -1/2$,
and the implicit constant is absolute. Since
$\zeta _{\fq (X)}={1\over (1-q^{-s})(1-q^{1-s})},$ we see that
$$\text{Res}_{s=1}\zeta _{\fq (X)}(s)={q\over (q-1)
\log (q)}.$$ 
Therefore
$$\aligned r&=\text{Res}_{s=1}f(s)\\
&={L_K(q^{-1})\over \zeta _K(2)}
\prod \Sb v\in\supp (\fc )\endSb (1+q^{-\deg (v)})^{-1}
\text{Res} _{s=1}\zeta _{\fq (X)}(s)\\
&=
{L_K(q^{-1})\over \zeta _K(2)}
\prod \Sb v\in\supp (\fc )\endSb (1+q^{-\deg (v)})^{-1}
{q\over (q-1)\log (q)}\\
&=
{J_Kq^{1-g_K}\over \zeta _K(2)(q-1)\log (q)}
\prod \Sb v\in\supp (\fc )\endSb (1+q^{-\deg (v)})^{-1}.\endaligned\tag 26$$
Using Lemma 0, we get
$$\aligned 
M&=\max \Sb \Re (s)=\epsilon -1/2\endSb \left\{\left | 
\prod _{v\not\in\supp (\fc )}(1-q^{-s\deg (v)})^{-1}(1-q^{-(s+1)\deg (v)})
\right |\right \}\\
&\ll
\max \Sb \Re (s)=\epsilon -1/2\endSb \left\{\left | 
\prod _{v\in\supp (\fc )}(1-q^{-s\deg (v)})(1-q^{-(s+1)\deg (v)})^{-1}\right |
\right \}\\
&\le
\prod _{v\in\supp (\fc )}(1+q^{-\epsilon +1/2)})
(1-q^{-\epsilon-1/2})^{-1}\\
&<
\prod _{v\in\supp (\fc )}(1+q^{1/2})(1-q^{-1/2})^{-1}\\
&\ll q^{\epsilon\deg (\fc )}.\endaligned\tag 27$$

We now set $j=m-2g_K+1$ in (25) and multiply through by $q^j.$ 
Using (26) and (27) gives
$$\sum\Sb \fd \ge 0\\ \deg (\fd )=m-2g_K+1\endSb \phi (\fd )\chic (\fd )
={J_Kq^{3-5g_K}\over \zetak (2)(q-1)}q^{2m}\prod\Sb v\in\supp (\fc )\endSb
(1+q^{-\deg (v)})^{-1}+O(q^{\epsilon\deg (\fc )}q^{\epsilon m}).$$
Using this estimate together with the Corollary to Proposition 3 and (0), we get
$$\aligned
\sum\Sb \fd\in\Div (K)\\ \fd \ge 0\\ \deg (\fd )=m-2g_K+1\endSb &
\sum\Sb [F:K]=2, \ q_F=q\\ \diff (F)=2\fd \endSb\chi (F/2\fc )\\
&=
\sum\Sb \fd\in\Div (K)\\ \fd \ge 0\\ \deg (\fd )=m-2g_K+1\endSb 
N(\fd )\chic (\fd )\\
&=
\sum\Sb \fd\in\Div (K)\\ \fd \ge 0\\ \deg (\fd )=m-2g_K+1\endSb 
2\phi (\fd )\chic (\fd )
+O\left (
\sum\Sb \fd\in\Div (K)\\ \fd \ge 0\\ \deg (\fd )=m-2g_K+1\endSb \chic (\fd )
\right )\\
&=
{2J_Kq^{3-5g_K}\over \zetak (2)(q-1)}q^{2m}\prod _{v\in\supp (\fc )}
(1+q^{-\deg (v)})^{-1}+O(q^{\epsilon m}q^{\epsilon\deg (\fc )}+q^m).
\endaligned$$

We remark that one can use a similar Tauberian argument for the case where
$q$ is odd (cf. [8, Proposition 17.2] which is the case where $\fc =0$),
though this leads to an estimate for the number of all square-free effective
divisors $\fd$ of fixed degree with $(\fc ,\fd )=0$. This gives Proposition 7
if the map $\psi$ in Proposition 2 has a trivial kernel, i.e., if $J_K$ is odd.
\enddemo

\head 5. Proofs of Theorems 1, 2 and 3\endhead
From the definitions in \S 1 we see that
$$\sum\Sb [F:K]=2\\ g_F=m,\ q_F=q_K\endSb L^*_F(q^{-s})L^*_F(q^{-t})=
\sum\Sb \fc \ge 0\endSb\sum\Sb \fa ,\fb\ge 0\\ \fa +\fb =\fc\endSb q^{-s\deg (\fa )}
q^{-t\deg (\fb )}\sum\Sb [F:K]=2\\ g_F=m, q_F=q_K\endSb \chi \big ( F/(\fa +\fb )\big )
,$$
with similar expressions including M\"obius functions for sums involving quotients
of $L$-functions. As indicated at the
end of \S 1 and demonstrated in Propositions 5-7,
the summands above where $\fc =\fa +\fb\in 2\Div (K)$ dominate.  We will
deal with these summands first.

\proclaim{Lemma 23} Suppose $K$ is a function field with field of
constants $\fq$ and $s,t\in\bc$ with $\Re (t),\Re (s)>1/2$. Then 
$$
\sum\Sb\fc\ge 0\endSb \sum\Sb \fa ,\fb\ge 0\\
\fa +\fb =2\fc\endSb q^{-s\deg (\fa )}q^{-t\deg (\fb )}
\prod _{v\in\supp (\fc )}(1+q^{-\deg (v)})^{-1}
=\zetak (2)\zetak (2s)\zetak (2t)\sigma _1(s,t),$$
$$
\sum\Sb\fc\ge 0\endSb \sum\Sb \fa ,\fb\ge 0\\
\fa +\fb =2\fc\endSb \mu (\fb ) q^{-s\deg (\fa )}q^{-t\deg (\fb )}
\prod _{v\in\supp (\fc )}(1+q^{-\deg (v)})^{-1}
=\zetak (2)\zetak (2s)\sigma _2(s,t),$$
and
$$
\sum\Sb\fc\ge 0\endSb \sum\Sb \fa ,\fb\ge 0\\
\fa +\fb =2\fc\endSb \mu (\fa )\mu (\fb ) q^{-s\deg (\fa )}q^{-t\deg (\fb )}
\prod _{v\in\supp (\fc )}(1+q^{-\deg (v)})^{-1}
=\zetak (2)\sigma _3(s,t).$$
\endproclaim

\demo{Proof} 
If $\fa+\fb =2\fc,$ then 
$-t\deg (\fb )-s\deg (\fa ) = -(t-s)\deg (\fb )-2s\deg (\fc )$.
Thus
$$\multline
\sum\Sb\fc\ge 0\endSb \sum\Sb \fa ,\fb\ge 0\\
\fa +\fb =2\fc\endSb \mu (\fb )q^{-s\deg (\fa )}q^{-t\deg (\fb )}
\prod _{v\in\supp (\fc )}(1+q^{-\deg (v)})^{-1}\\
=\sum\Sb\fc\ge 0\endSb \sum\Sb 0\le\fb\le 2\fc\endSb
\mu (\fb )q^{-(t-s)\deg (\fb )}q^{-2s\deg (\fc )}
\prod _{v\in\supp (\fc )}(1+q^{-\deg (v)})^{-1}\endmultline\tag 28$$
and similarly for sums where the $\mu (\fb )$ term is absent.
If $\fa +\fb =2\fc$ and $\mu (\fa )\mu (\fb )\neq 0$, then $\fa =\fb =\fc$
and $\mu (\fa )\mu (\fb )=|\mu (\fc )|$. Set
$$ \aligned
\theta _1(\fc )&=q^{-2s\deg (\fc )}\prod \Sb v\in\supp (\fc )\endSb
(1+q^{-\deg (v)})^{-1}
\sum\Sb 0\le \fb\le 2\fc\endSb q^{-(t-s)\deg (\fb )}\\
\theta _2(\fc )&=q^{-2s\deg (\fc )}\prod \Sb v\in\supp (\fc )\endSb
(1+q^{-\deg (v)})^{-1}
\sum\Sb 0\le \fb\le 2\fc\endSb \mu (\fb )q^{-(t-s)\deg (\fb )}\\
\theta _3(\fc )&=|\mu (\fc )|q^{-(s+t)\deg (\fc )}\prod\Sb v\in\supp (\fc )\endSb
(1+q^{-\deg (v)})^{-1}.\endaligned$$
Note that $\theta _i (\fc +\fd )=\theta _i(\fc )\theta  _i(\fd )$
whenever $(\fc ,\fd )=0$ for all $i$. Thus 
$$\sum \Sb\fc\ge 0\endSb\theta _i(\fc )=\prod \Sb v\in M(K)\endSb \left (
1+\sum _{k=1}^{\infty}\theta _i (kv)\right ).\tag 29$$

Suppose that $s\neq t$. Then for all positive integers $k$ and all places
$v$
$$\aligned
\theta _1(kv)&=(1+q^{-\deg (v)})^{-1}q^{-2ks\deg (v)}\sum _{i=0}^{2k} q^{-(t-s)i
\deg (v)}\\
&=(1+q^{-\deg (v)})^{-1}q^{-2ks\deg (v)}(1-q^{-(t-s)(2k+1)\deg (v)})(1-q^{-(t-s)
\deg (v)})^{-1}\\
&=(1+q^{-\deg (v)})^{-1}(q^{-(2k+1)s\deg (v)}-q^{-(2k+1)t\deg (v)})(q^{-s\deg (v)}-
q^{-t\deg (v)})^{-1}.\endaligned$$
Therefore
$$\aligned\sum _{k=1}^{\infty}\theta _1(kv)&=
(1+q^{-\deg (v)})^{-1}(q^{-s\deg (v)}-q^{-t\deg (v)})^{-1}\\
&\qquad\times \Big (
q^{-3s\deg (v)}(1-q^{-2s\deg (v)})^{-1}-q^{-3t\deg (v)}(1-q^{-2t\deg (v)})^{-1}\Big )\\
&=
(1+q^{-\deg (v)})^{-1}(1-q^{-2s\deg (v)})^{-1}(1-q^{-2t\deg (v)})^{-1}
(q^{-s\deg (v)}-q^{-t\deg (v)})\\
&\qquad\times\Big ( q^{-3s\deg (v)}-q^{-3t\deg (v)}-
q^{-2(s+t)\deg (v)}(q^{-s\deg (v)}-q^{-t\deg (v)}\Big )\\
&=
(1-q^{-2\deg (v)})^{-1}(1-q^{-2s\deg (v)})^{-1}(1-q^{-2t\deg (v)})^{-1}(1-q^{-\deg (v)})
\\
&\qquad\times
\Big (q^{-2s\deg (v)}+q^{-(s+t)\deg (v)}+q^{-2t\deg (v)}-q^{-2(s+t)\deg (v)}\Big ),
\endaligned$$
so that by (29)
$$\aligned
\sum\Sb \fc\ge 0\endSb \theta _1(\fc )&=
\zetak (2)\zetak (2s)\zetak (2t)\\
&\times
\prod _{v\in M(K)}\Big ( (1-q^{-2\deg (v)})(1-q^{-2s\deg (v)})(1-q^{-2t\deg (v)})\\
&\qquad
+(1-q^{-\deg (v)})(q^{-2s\deg (v)}+q^{-(s+t)\deg (v)}+q^{-2t\deg (v)}-q^{-2(s+t)
\deg (v)})\Big )\\
&=\zetak (2)\zetak (2s)\zetak (2t)\sigma _1(s,t).\endaligned\tag 30$$

Now if $s=t$, then via term-by-term differentiation of
$$\sum _{k=1}^{\infty}x^{2k+1}=x\sum _{k=1}^{\infty}x^{2k}=x^3\sum _{k=0}
^{\infty} x^{2k}= {x^3\over 1-x^2}$$
shows  that
$$\aligned \sum _{k=1}^{\infty} \theta _1(kv)&=(1+q^{-\deg (v)})^{-1}
\sum _{k=1}^{\infty} (2k+1)q^{-2sk\deg (v)}\\
&= (1+q^{-\deg (v)})^{-1} {3q^{-2s\deg (v)}-q^{-4s\deg (v)}\over
(1-q^{-2s\deg (v)})^2}\\
&=(1-q^{-2\deg (v)})^{-1}(1-q^{-2s\deg (v)})^{-2}(1-q^{-\deg (v)})(3q^{-2s\deg (v)}-
q^{-4s\deg (v)}).\endaligned$$
Using this together with (29) yields
$$\aligned
\sum \Sb \fc\ge 0\endSb \theta _1(\fc )
&=\zetak (2)\zetak (2s)^2\\
&\times
\prod\Sb v\in M(K)\endSb (1-q^{-2\deg (v)})(1-q^{-2s\deg (v)})^2+
(1-q^{-\deg (v)})(3q^{-2s\deg (v)}-q^{-4s\deg (v)})\\
&=\zetak (2)\zetak (2s)\zetak (2t)\sigma _1(s,t)\endaligned\tag 30'$$
when $s=t$. Together,  (28) and (30) (when $s\neq t$), or (28) and  (30') 
(when $s=t$) proves the first part of the lemma.

For all positive integers $k$ and $v\in M(K)$ we have
$$\theta _2 (kv)=
q^{-2sk\deg (v)}(1+q^{-\deg (v)})^{-1}(1-q^{-(t-s)\deg (v)}),$$
so that
$$\aligned\sum _{k=1}^{\infty}\theta (kv)&=(1+q^{-\deg (v)})^{-1}(1-q^{-(t-s)\deg (v)})
q^{-2s\deg (v)}(1-q^{-2s\deg (v)})^{-1}\\
&=(1-q^{-2\deg (v)})^{-1}(1-q^{-2s\deg (v)})^{-1}(1-q^{-\deg (v)})(q^{-2s\deg (v)}-
q^{-(t+s)\deg (v)}).\endaligned$$
Whence by (29)
$$\aligned \sum\Sb \fc \ge 0\endSb \theta _2(\fc )&=\zetak (2)\zetak (2s)\\
&\times
\prod _{v\in M(K)}\Big ( (1-q^{-2\deg (v)})(1-q^{-2s\deg (v)})\\
&\qquad + (1-q^{-\deg (v)})(q^{-2s\deg (v)}-q^{-(t+s)\deg (v)})\Big )\\
&=\zetak (2)\zetak (2s)\sigma _2(s,t).\endaligned\tag 31$$
The second part of Lemma 23 follows from (28) and (31).

Finally, we easily see that for all positive integers $k$ and places $v$
$$\theta _3(kv)=\cases q^{-(s+t)\deg (v)}(1+q^{-\deg (v)})^{-1}&\text{if $k=1$,}\\
0&\text{if $k>1$.}\endcases$$
Thus
$$\aligned \prod _{v\in M(K)}\left ( 1+\sum _{k=1}^{\infty}\theta _3(kv)\right )&=
\prod _{v\in M(K)}
1+q^{-(s+t)\deg (v)}(1+q^{-\deg (v)})^{-1}\\
&=\prod _{v\in M(K)}
(1+q^{-\deg (v)})^{-1}( 1+q^{-\deg (v)}+q^{-(s+t)\deg (v)} )\\
&=\prod _{v\in M(K)}
(1-q^{-2\deg (v)})^{-1}(1-q^{-\deg (v)})(1+q^{-\deg (v)}+q^{-(s+t)\deg (v)})\\
&=\zetak (2)\sigma _3(s,t).\endaligned$$
The last part of the lemma follows from this and (29).
\enddemo

\proclaim{Lemma 24} Let $K$ be a function field with field of constants
$\fq$. Suppose $m$ is a positive integer and $\epsilon >0$. Then for all
$s,t\in\bc$ with $\Re (t),\Re (s)>1/2+\epsilon$
$$\multline
\sum \Sb \fc\ge 0\endSb \sum\Sb \fa ,\fb\in\Div (K)\\
\fa ,\fb \ge 0\\ \fa +\fb =2\fc\endSb q^{-s\deg (\fa )}q^{-t\deg (\fb )}
\sum\Sb [F:K]=2\\ g_F=m,\ q_F=q\endSb \chi (F/2\fc )\\
={2J_Kq^{3-5g_K}\zetak (2s)\zetak (2t)\over q-1}\sigma _1(s,t)q^{2m}
+\cases O(q^{(1/2+\epsilon )m})&\text{if $q$ is odd,}\\
O(q^m)&\text{if $q$ is even,}\endcases\endmultline$$
$$\multline
\sum \Sb \fc\ge 0\endSb \sum\Sb \fa ,\fb\in\Div (K)\\
\fa ,\fb \ge 0\\ \fa +\fb =2\fc\endSb \mu (\fb )q^{-s\deg (\fa )}
q^{-t\deg (\fb )}\sum\Sb [F:K]=2\\ g_F=m,\  q_F=q\endSb \chi (F/2\fc )\\
={2J_Kq^{3-5g_K}\zetak (2s)\over q-1}\sigma _2(s,t)q^{2m}
+\cases O(q^{(1/2+\epsilon )m})&\text{if $q$ is odd,}\\
O(q^m)&\text{if $q$ is even,}\endcases\endmultline$$
and
$$\multline
\sum \Sb \fc\ge 0\endSb \sum\Sb \fa ,\fb\in\Div (K)\\
\fa ,\fb \ge 0\\ \fa +\fb =2\fc\endSb\mu (\fa )\mu (\fb )q^{-s\deg (\fa )}
q^{-t\deg (\fb )}\sum\Sb [F:K]=2\\ g_F=m,\  q_F=q\endSb \chi (F/2\fc )\\
={2J_Kq^{3-5g_K}\over q-1}\sigma _3(s,t)q^{2m}
+\cases O(q^{(1/2+\epsilon )m})&\text{if $q$ is odd,}\\
O(q^m)&\text{if $q$ is even,}\endcases\endmultline$$
where the implicit constants depend only on $K$ and $\epsilon$.
\endproclaim

\demo{Proof} The main terms follow directly from Proposition 7 and Lemma 23.
As for the sums over the error terms in Proposition 7, let ${\frak m}=
\min\{\Re (s),\Re (t)\}$. Then for all $\delta >0$ we have
$$\aligned\sum\Sb\fc\ge 0\endSb
\sum\Sb \fa ,\fb\ge 0\\ 
\fa +\fb=2\fc \endSb q^{-\Re (s)\deg (\fa )}q^{-\Re (t)\deg (\fb )}&\le
\sum \Sb \fc\ge 0\endSb q^{-2{\frak m}\deg (\fc )}\sum\Sb 0\le\fd\le 2\fc\endSb 
1\\
&\ll \sum \Sb \fc\ge 0\endSb q^{-(2{\frak m}-\delta /2 )\deg (\fc )},\endaligned
\tag 32$$
by Lemma 0, where the implicit constant depends only on $K$ and $\delta$.
Now set $2\delta = {\frak m}-1/2$. Then by (32) and replacing the $\epsilon$ in
Proposition 7 with $\delta /2$, we have
$$\aligned
\sum\Sb\fc\ge 0\endSb \sum\Sb \fa ,\fb \ge 0\\ \fa +\fb =2\fc\endSb
q^{-\Re (s)\deg (\fa )}q^{-\Re (t)\deg (\fb )}q^{(\delta /2)\deg (\fc )}
&\ll
\sum\Sb \fc\ge 0\endSb q^{-(2{\frak m}-\delta )\deg (\fc )}\\
&\ll 1.\endaligned\tag 33$$
The error terms in Lemma 24 follow from Proposition 7 and (33).
\enddemo

With the sums over the main terms done, we now turn to the sums over the
error terms, i.e., the sums
where $\fa +\fb\not\in 2\Div (K)$. When no M\"obius function appears with
a given divisor, then we can use the following estimate when the degree of
that divisor is relatively large in terms of the parameter $m$.

\proclaim{Lemma 25} Suppose $K$ is a function field with field
of constants $\fq$ and $m$ is a non-negative integer. 
Let $\epsilon >0$ and suppose
$s,t\in\bc$ with $\Re (t), \Re (s)> 1/2+\epsilon$. Then for any
quadratic extension $F\supset K$ with $g_F=m$ and $q_F=q$
$$\sum\Sb \fb\ge 0\endSb q^{-\Re (t)\deg (\fb )}
\left |\sum\Sb \fa \ge 0\fa +\fb\not\in 2\Div (K)\\ \deg (\fa )> 2m-2g_K\endSb 
q^{-s\deg (\fa )}\chi \big ( F/\fa )\big )\right |
\ll q^{-2m(\Re (s)- 1/2)},$$
where the implicit constant depends only on $K$ and $\epsilon$. 
\endproclaim

\demo{Proof}
Fix a quadratic extension $F$ of $K$ with $g_F=m$
and $q_F=q$.  
For the moment, fix a divisor $\fb \ge 0$ and write $\fb =\fb '+2\fb ''$ where
$\fb ',\fb ''\ge 0$ and $\fb '$ is square-free. Let $n>2m-2g_K$ and suppose that
$\fa$ is an effective divisor such that $\deg (\fa )=n$ and
$\fa +\fb =2\fc$ for some divisor $\fc$. Then we must have
$\fa =\fb ' +2\fc '$ for some effective divisor $\fc '$ with
$2\deg (\fc ')=n-\deg (\fb ')$. By a theorem of Weil (see [8, Theorem 9.16B]),
$L_F^*(q^{-s})$ is a polynomial of degree $2g_F-2g_K$ in $q^{-s}$. 
Since $n>2m-2g_K=2g_F-2g_K$, 
$$\aligned 0&=\sum\Sb\fa\ge 0\\ \deg (\fa )=n\endSb \chi (F/\fa )\\
&=\sum\Sb\fa\ge 0\\ \deg (\fa )=n\\ \fa =\fb '+2\fc '\endSb\chi (F/\fa )
+\sum\Sb\fa\ge 0\\ \deg (\fa )=n\\ \fa +\fb\not\in 2\Div (K)\endSb
\chi (F/\fa )\\
&=\chi (F/\fb ')\sum\Sb \fc '\ge 0\\ 2\deg (\fc ')=n-\deg (\fb ')\endSb\chi (F/2\fc ')
+\sum\Sb\fa\ge 0\\ \deg (\fa )=n\\ \fa +\fb\not\in 2\Div (K)\endSb
\chi (F/\fa ),\endaligned$$
so that by (0)
$$\aligned \left | \sum\Sb \fa\ge 0\\ \deg (\fa )=n\\ \fa +\fb\not\in
2\Div (K)\endSb \chi (F/\fa )\right |&=
\left | \chi (F/\fb ')\sum\Sb \fc '\ge 0\\ 2\deg (\fc ')=n-\deg (\fb ')\endSb
\chi (F/2\fc ')\right |\\
&\le\sum\Sb \fc '\ge 0\\ 2\deg (\fc ')=n-\deg (\fb ')\endSb 1\\
&\ll q^{(n-\deg (\fb '))/2},\endaligned$$
where the implicit constant depends on $K$ only.
Using this, we see that 
$$
\left |\sum\Sb\fa\ge 0\\ \deg (\fa )=n\\ \fa +\fb\not\in 2\Div (K)\endSb
\chi (F/\fa)q^{-s\deg (\fa )}\right |
\ll q^{-\deg (\fb ')/2}q^{-n(\Re (s)-1/2)}.\tag 34$$
Since $\Re (s),\ \Re (t)> 1/2+\epsilon,$
$$\sum\Sb n>2m-2g_K\endSb q^{-n(\Re (s)-1/2)}\ll q^{-2m(\Re (s)-1/2)}\tag 35$$
and
$$\aligned
\sum\Sb\fb\ge 0\endSb q^{-\Re (t)\deg (\fb )}q^{-\deg (\fb ')/2}&=
\sum\Sb \fb\ge 0\endSb q^{-(\Re (t)+1/2)\deg (\fb ')}q^{-2\Re (t)\deg (\fb '')}\\
&\le\sum\Sb\fb '\ge 0\endSb q^{-(\Re (t)+1/2)\deg (\fb ')}\sum\Sb\fb ''\ge 0\endSb
q^{-2\Re (t)\deg (\fb '')}\\
&\ll 1.\endaligned\tag 36$$
The lemma follows from (34)-(36).\enddemo

\demo{Proof of Theorem 1} By Proposition 7 (specifically, the estimate for the
number of quadratic extensions of genus $m$), Lemma 25 and the multiplicative
nature of $\chi (F/*)$, 
$$\multline
\sum\Sb [F:K]=2\\ g_F=m,\ q_F=q\endSb \left |\sum\Sb \fa ,\fb \ge 0\\ 
\fa +\fb\not\in 2\Div (K)\\ \deg (\fa )> 2m-2g_K \ \text{or}\ \deg (\fb )>2m-2g_K
\endSb q^{-s\deg (\fa )}q^{-t\deg (\fb )}\chi\big ( F/(\fa +\fb )\big )\right |\\
\ll q^{2m(3/2-\Re (s))}+q^{2m(3/2-\Re (t))}.\endmultline\tag 37$$

Suppose that $q$ is odd. By Proposition 5 with $\fa +\fb$ in place of $\fc$, 
$$\left |\sum\Sb [F:K]=2\\ g_F=m, q_F=q\endSb \chi \big ( F/(\fa +\fb )\big )\right |
\ll q^{m+(\epsilon +1/4)\deg (\fa +\fb )}.$$
This in conjunction with (0) yields
$$\multline
\left |\sum\Sb \fa ,\fb \ge 0\\ \fa +\fb\not\in 2\Div (K)\\ \deg (\fa ),\deg (\fb )
\le 2m-2g_K\endSb q^{-s\deg (\fa )}q^{-t\deg (\fb )}\sum\Sb [F:K]=2\\ g_F=m,\ q_F=q
\endSb \chi \big ( F/(\fa +\fb )\big )\right |\\
\ll q^{m}
\sum\Sb \fa\ge 0\\ \deg (\fa )\le 2m-2g_K\endSb q^{(\epsilon +1/4-\Re (s))\deg (\fa )}
\sum\Sb \fb\ge 0\\ \deg (\fb )\le 2m-2g_K\endSb q^{(\epsilon +1/4-\Re (t))\deg (\fb )}
\\
\ll q^m(1+q^{2m(5/4 +\epsilon -\Re (s))})(1+q^{2m(5/4 +\epsilon -\Re (t))})
.\endmultline\tag 38$$

Suppose that $q$ is even. By Proposition 6 with $\fa +\fb$ in place of $\fc$,
$$\left |\sum\Sb [F:K]=2\\ g_F=m, q_F=q\endSb \chi \big ( F/(\fa +\fb )\big )\right |
\ll q^{m+\epsilon\deg (\fa +\fb )}.$$
This in conjunction with (0) yields
$$\multline
\left |\sum\Sb \fa ,\fb \ge 0\\ \fa +\fb\not\in 2\Div (K)\\ \deg (\fa ),\deg (\fb )
\le 2m-2g_K\endSb q^{-s\deg (\fa )}q^{-t\deg (\fb )}\sum\Sb [F:K]=2\\ g_F=m,\ q_F=q
\endSb \chi \big ( F/(\fa +\fb )\big )\right |\\
\ll q^{m}
\sum\Sb \fa\ge 0\\ \deg (\fa )\le 2m-2g_K\endSb q^{(\epsilon -\Re (s))\deg (\fa )}
\sum\Sb \fb\ge 0\\ \deg (\fb )\le 2m-2g_K\endSb q^{(\epsilon -\Re (t))\deg (\fb )}
\\
\ll q^m(1+q^{2m(1 +\epsilon -\Re (s))})(1+q^{2m(1 +\epsilon -\Re (t))})
.\endmultline\tag 39$$

Theorem 1 follows from Lemma 24 and (37)-(39).
\enddemo

\demo{Proof of Theorem 2}
Since we have a $\mu (\fb )$ factor now, instead of (37) we have
$$\multline
\sum\Sb [F:K]=2\\ g_F=m,\ q_F=q\endSb \left |\sum\Sb \fa ,\fb \ge 0\\ 
\fa +\fb\not\in 2\Div (K)\\ \deg (\fa )> 2m-2g_K\endSb\mu (\fb )
q^{-s\deg (\fa )}q^{-t\deg (\fb )}\chi\big ( F/(\fa +\fb )\big )\right |\\
\ll q^{2m(3/2-\Re (s))}.\endmultline\tag 40$$

Suppose $q$ is odd. Using Proposition 5 with $\fa +\fb$ in place of $\fc$
we have
$$\multline
\left |\sum\Sb \fa ,\fb \ge 0\\ \fa +\fb\not\in 2\Div (K)\\ \deg (\fa ) +\deg (\fb )
\le 4m\\ \deg (\fa )\le 2m-2g_K\endSb \mu (\fb )
q^{-s\deg (\fa )}q^{-t\deg (\fb )}\sum\Sb [F:K]=2\\ g_F=m,\ q_F=q
\endSb \chi \big ( F/(\fa +\fb )\big )\right |\\
\ll q^{m}
\sum\Sb \fa\ge 0\\ \deg (\fa )\le 2m-2g_K\endSb q^{(\epsilon +1/4-\Re (s))\deg (\fa )}
\sum\Sb \fb\ge 0\\ \deg (\fb )\le 4m-\deg (\fa )\endSb 
q^{(\epsilon +1/4-\Re (t))\deg (\fb )}\\
\ll q^m
\sum\Sb \fa\ge 0\\ \deg (\fa )\le 2m-2g_K\endSb q^{(\epsilon +1/4-\Re (s))\deg (\fa )}
(1+q^{(5/4+\epsilon -\Re (t))(4m-\deg (\fa ))})\\
\ll q^m\big (1+q^{2m(5/4+\epsilon -\Re (s))}+q^{2m(5/2+2\epsilon -2\Re (t))}+
q^{2m(5/2+2\epsilon -\Re (s)-\Re (t))}\big ) ,\endmultline\tag 41$$
and similarly, since $\Re (t)>1+\epsilon$
$$\multline
\left |\sum\Sb \fa ,\fb \ge 0\\ \fa +\fb\not\in 2\Div (K)\\ \deg (\fa )+\deg (\fb )
> 4m\\ \deg (\fa )\le 2m-2g_K\endSb \mu (\fb )
q^{-s\deg (\fa )}q^{-t\deg (\fb )}\sum\Sb [F:K]=2\\ g_F=m,\ q_F=q
\endSb \chi \big ( F/(\fa +\fb )\big )\right |\\
\ll q^{2m}
\sum\Sb \fa\ge 0\\ \deg (\fa )\le 2m-2g_K\endSb q^{-\Re (s)\deg (\fa )}
\sum\Sb \fb\ge 0\\ \deg (\fb )> 4m-\deg (\fa )\endSb q^{-\Re (t)\deg (\fb )}\\
\ll q^{2m}q^{-4m(\Re (t)-1)}
\sum\Sb \fa\ge 0\\ \deg (\fa )\le 2m-2g_K\endSb q^{(\Re (t)-1-\Re (s))\deg (\fa )}\\
\ll q^m\big (q^{2m(5/2+2\epsilon -2\Re (t))}+q^{2m(5/2-\Re (s)-\Re (t))}\big ).
\endmultline\tag 42$$

Suppose that $q$ is even. By Proposition 6 with $\fa +\fb$ in place of $\fc$
and $\epsilon /2$ in place of $\epsilon$, and using $\Re (t)> 1+\epsilon$ we get
$$\multline
\left |\sum\Sb \fa ,\fb \ge 0\\ \fa +\fb\not\in 2\Div (K)\\ \deg (\fa )
\le 2m-2g_K\endSb \mu (\fb )q^{-s\deg (\fa )}q^{-t\deg (\fb )}
\sum\Sb [F:K]=2\\ g_F=m,\ q_F=q
\endSb \chi \big ( F/(\fa +\fb )\big )\right |\\
\ll q^{m}
\sum\Sb \fa\ge 0\\ \deg (\fa )\le 2m-2g_K\endSb q^{(\epsilon -\Re (s))\deg (\fa )}
\sum\Sb \fb\ge 0\endSb q^{(\epsilon -\Re (t))\deg (\fb )}\\
\ll q^m(1+q^{2m(1 +\epsilon -\Re (s))})
.\endmultline\tag 43$$

Theorem 2 follows from Lemma 24 and (40)-(43).
\enddemo

\demo{Proof of Theorem 3}
Since we now have both a $\mu (\fa )$ and $\mu (\fb )$ term, we are unable to
utilize Lemma 25. 

Suppose $q$ is odd. Then by Proposition 5  we have
$$\multline
\left |\sum\Sb \fa ,\fb \ge 0\\ \fa +\fb\not\in 2\Div (K)\\ \deg (\fa )+\deg (\fb )
\le 4m\endSb
\mu (\fa )\mu (\fb )q^{-s\deg (\fa )}q^{-t\deg (\fb )}\sum\Sb [F:K]=2\\ g_F=m,\ q_F=q
\endSb \chi \big ( F/(\fa +\fb )\big )\right |\\
\ll q^{m}
\sum\Sb \fa\ge 0\\ \deg (\fa )\le 4m\endSb q^{(\epsilon  +1/4-\Re (s))\deg (\fa )}
\sum\Sb \fb\ge 0\\ \deg (\fb )\le 4m-\deg (\fa )
\endSb q^{(\epsilon +1/4-\Re (t))\deg (\fb )}\\
\ll q^m
\sum\Sb \fa\ge 0\\ \deg (\fa )\le 4m\endSb q^{(\epsilon /2 +1/4-\Re (s))\deg (\fa )}
\big ( 1+q^{(4m -\deg (\fa )(5/4+\epsilon -\Re (t))}\big )\\
\ll q^m\big (1+q^{2m(5/2+\epsilon -2\Re (s))}+q^{2m(5/2+\epsilon -2\Re (t))}\big )
.\endmultline\tag 44$$
Similarly, we get
$$\multline
\left |\sum\Sb \fa ,\fb \ge 0\\ \fa +\fb\not\in 2\Div (K)\\ \deg (\fa )+\deg (\fb )
>4m\endSb
\mu (\fa )\mu (\fb )q^{-s\deg (\fa )}q^{-t\deg (\fb )}\sum\Sb [F:K]=2\\ g_F=m,\ q_F=q
\endSb \chi \big ( F/(\fa +\fb )\big )\right |\\
\ll q^m\big (1+q^{2m(5/2-2\Re (s))}+q^{2m(5/2-2\Re (t))}\big )
.\endmultline\tag 45$$

Suppose $q$ is even. Then by Proposition 6 with $\epsilon$ replaced by 
$\epsilon /2$ and using $\Re (s),\ \Re (t)> 1+\epsilon$, we have
$$\multline
\left |\sum\Sb \fa ,\fb \ge 0\\ \fa +\fb\not\in 2\Div (K)\endSb
\mu (\fa )\mu (\fb )q^{-s\deg (\fa )}q^{-t\deg (\fb )}\sum\Sb [F:K]=2\\ g_F=m,\ q_F=q
\endSb \chi \big ( F/(\fa +\fb )\big )\right |\\
\ll q^m
\sum\Sb \fa\ge 0\endSb q^{(\epsilon /2  -\Re (s))\deg (\fa )}
\sum\Sb \fb\ge 0\endSb q^{(\epsilon /2  -\Re (t))\deg (\fb )}\\
\ll q^m.\endmultline\tag 46$$

Theorem 3 follows from Lemma 24 and (44)-(46).
\enddemo

\enddocument

\enddocument